\documentclass{amsart}  
\usepackage{amsmath,amsfonts,mathrsfs,enumerate}

\usepackage{graphicx}

\newcommand{\sA}{{\mathcal A}}
\newcommand{\sB}{{\mathcal B}}
\newcommand{\bsC}{{\boldsymbol{\mathcal{C}}}} 
\newcommand{\sC}{{\mathcal{C}}}
\newcommand{\sD}{{\mathcal D}}
\newcommand{\sE}{{\mathcal E}}
\newcommand{\sF}{{\mathcal F}}

\newcommand{\sI}{{\mathcal I}}
\newcommand{\sM}{{\mathcal M}}
\newcommand{\sN}{{\mathcal N}}
\newcommand{\sY}{{\mathcal Y}}

\newcommand{\vh}{{\mathbf{h}}}
\newcommand{\hz}{{\hat{z}}}
\newcommand{\cz}{{\overset{\scriptscriptstyle\circ}{z}}}
\newcommand{\tg}{{\tilde{g}}}
\newcommand{\tb}{{\tilde{b}}}
\newcommand{\tfw}{{\tilde{f_1}}}
\newcommand{\tfp}{{\tilde{f_2}}}

\newcommand{\myspace}{\qquad\qquad\qquad}

\newtheorem{theorem}{Theorem}[section]
\newtheorem{lemma}[theorem]{Lemma}
\newtheorem{proposition}[theorem]{Proposition}
\newtheorem{remark}{Remark}

\newtheorem{assumption}[theorem]{Assumption}

\newtheorem{definition}[theorem]{Definition}

\numberwithin{equation}{section}

\newcommand{\scrA}{\mathscr{A}}
\newcommand{\grad}{\nabla}

\newcommand{\bbR}{\mathbb{R}}
\newcommand{\bbN}{\mathbb{N}}
\newcommand{\p}{\partial}

\newcommand{\al}{\alpha}
\newcommand{\bet}{\beta}
\newcommand{\Del}{\Delta}
\newcommand{\del}{\delta}
\newcommand{\eps}{\epsilon}
\newcommand{\beps}{\bar{\eps}}

\newcommand{\Om}{\Omega}
\newcommand{\Gam}{\Gamma}
\newcommand{\gam}{\gamma}
\newcommand{\kap}{\kappa}
\newcommand{\lam}{\lambda}

\newcommand{\sig}{\sigma}
\newcommand{\zet}{\zeta}
\renewcommand{\th}{\theta}

\newcommand{\intQT}{\int_{Q_T}}
\newcommand{\intOm}{\int_\Omega}
\newcommand{\intT}{\int_0^T}
\newcommand{\intst}{\int_s^t}
\newcommand{\intsT}{\int_s^{s+T}}

\newcommand{\Dn}[1]{\frac{\partial #1}{\partial \nu}}
\renewcommand{\lg}{\langle}
\newcommand{\rg}{\rangle}
\newcommand{\into}{\hookrightarrow}

\newcommand{\zmn}{z^{m,n}}
\newcommand{\vmn}{v^{m,n}}
\newcommand{\dfn}{:=}
\renewcommand{\div}{\operatorname{div}}
\newcommand{\barT}{\overline{T}}

\title[Attractors for a coupled wave/plate system]
{Finite dimensional attractor for a composite system of wave/plate 
equations with localised damping}

\author[F. Bucci, D. Toundykov]{}

\thanks{This research was partly performed while both authors were 
visiting the Centro di Ricerca Matematica Ennio De Giorgi of the
Scuola Normale Superiore in Pisa, within the program {\em Research in pairs}.
The support (and kindness) of the Director Mariano Giaquinta, as well as the 
co-operative assistance of the secretarial staff of the Centro are warmly acknowledged.
\\
\indent
The research of F.~Bucci was partially supported by the Universit\`a degli Studi
di Firenze, 
and by the Italian MIUR, within the project 2007WECYEA 
(``Metodi di viscosit\`a, metrici e di teoria del controllo in equazioni alle derivate parziali non lineari'').
\\
\indent
The research of D.~Toundykov was partially supported by the National Science Foundation under Grant DMS-0908270.}

\begin{document}

\maketitle

\centerline{\scshape Francesca Bucci}
\medskip
{\footnotesize
\centerline{Universit\`a degli Studi di Firenze}
\centerline{Dipartimento di Matematica Applicata}
\centerline{Firenze, 50139, Italy}
} 

\bigskip

\centerline{\scshape Daniel Toundykov}
\medskip
{\footnotesize
\centerline{University of Nebraska--Lincoln}
\centerline{Department of Mathematics}
\centerline{Lincoln, NE 68588, USA}
}

\medskip

\date{}

\begin{quote}{\normalfont\fontsize{8}{10}\selectfont
{\bfseries Abstract.}
The long-term behaviour of solutions to a model for acoustic-structure interactions is addressed; the system is comprised of coupled semilinear wave (3D) and plate equations 
with nonlinear damping and critical sources. The questions of interest are: existence of 
a global attractor for the dynamics generated by this composite system, as well as dimensionality and regularity of the attractor.
A distinct and challenging feature of the problem is the geometrically restricted dissipation on the wave component of the system. It is shown that the existence of a global attractor of finite fractal dimension---established in a previous work by Bucci, Chueshov and Lasiecka (Comm. Pure Appl. Anal., 2007) only in the presence of \emph{full interior} acoustic damping---holds even in the case of \emph{localised} dissipation. This nontrivial generalization is inspired by and consistent with the recent advances in the study of wave equations with nonlinear localised damping.
\par}
\end{quote}

\section{Introduction}
This paper focuses on the long-term behaviour of a system of partial differential equations (PDE's, or simply PDE) modeling acoustic waves in a three-dimensional domain $\Omega$, and their interaction with the elastic vibrations of a part $\Gamma_0 \subset \Gam:=\p\Om$ of  the boundary.
The corresponding PDE problem (\eqref{e:pde-system-0} below) is comprised of semilinear wave and plate equations on domains  $\Omega$ and $\Gam_0$ respectively; a detailed description of the system will be given in the next section.

PDE models for acoustic-structure interactions have received a great deal of attention in the past decade, mainly in the context of control theory, in connection with (but not limited to) the problem of reducing the noise inside an acoustic chamber and stabilizing its vibrating walls.
As in the case of a single PDE, the sought-after properties---on a theoretical level,
yet with implications for the computational methods---are (besides well-posedness): stability, controllability, and solvability of the associated optimal control problems.
However, the presence of two (or more) coupled equations, of different type (usually, hyperbolic and parabolic) and/or acting on manifolds of different dimensions, renders the PDE analysis of the dynamics far more complex.
\\
It is beyond this work's scope to provide a comprehensive account of the literature 
on control problems for such interconnected PDE systems. However, for completeness and 
the reader's convenience some notable contributions are listed below.

Deep insight into the physical origin of PDE models for acoustic-structure interactions
is provided by \cite{morse-ingard}.
A more recent reference on analytical methods for modeling acoustic problems 
is \cite{howe}.

A general reference for the mathematical control theory of established 
coupled PDE systems for acoustic-structure interactions is the treatise 
\cite{las-cbms}, which also contains an extensive bibliography, from which 
we explicitly cite \cite{avalos-las-1} (PDE analysis and optimal control), \cite{avalos-las-2,avalos-las-3} (controllability), and \cite{cathe,las-cathe} (uniform stability), besides the former \cite{dimitriadis-etal}, \cite{beale}, \cite{littman-liu}.

More recent advances include, for instance (without any claim of completeness), 
\cite{las-trig-se-1}, \cite{las-tuff-2} and \cite{acquistapace-etal} (general theory of quadratic optimal control and of Riccati equations), 
\cite{bucci-1,bucci-2} (sharp interior/boundary regularity, with application to optimal control), \cite{cagnol-etal} (stabilization). 
For a survey of results on exact boundary controllability and uniform boundary stabilizability by differential geometric methods, see \cite{gulliver-etal}.
The stability of an interesting 2D structural-acoustic model was discussed
in \cite{grobbelaar, grobbelaar:2}. 

\smallskip
The focus of this study is on the long-term behavior of the nonlinear dynamics
generated by the system \eqref{e:pde-system-0} below, such as existence of a global attractor and its properties: geometry, fractal dimension, regularity.
It was shown in \cite{bu-chu-las} that in the presence of an acoustic dissipation distributed over the entire domain $\Omega$, existence of a global attractor 
is guaranteed under some basic assumptions on all the nonlinearities, including the most relevant case when the function modeling the acoustic source has a {\em critical exponent}. Furthermore, even under a critical, i.e.~non-compact, perturbation the attractor has a finite fractal dimension, provided the dissipation terms acting on either equation satisfy suitable conditions; in particular, the acoustic dissipation must be {\em linearly bounded}.

In light of the recent advances in the study of the long-term behaviour of wave equations with nonlinear damping \cite{chu-las-tou, chu-las-tou2}, the present work aims at showing that dissipation active only in a neighbourhood of $\Gamma_0$---combined with a nonlinear interior plate damping, as in the PDE model addressed in \cite{bu-chu-las}---still guarantees the existence of a finite-dimensional global attractor.

With this as background motivation, let us now discuss the major challenges which arise in the present context of 
\begin{itemize}
	\item  coupled hyperbolic PDE's (wave and plate equations),
	
	\item  \emph{localised damping in conjunction with a critical exponent-source} on the wave equation,
	
	\item  a critical-exponent source on the plate equation,
\end{itemize}
and then describe the principal theoretical and technical arguments utilized to overcome them.

\subsection{New challenges}\label{ss:newchallenges}
Critical sources alone have been well-known to present an intrinsic obstacle to formation 
of attractors. At the level of critical-exponents the compactness of corresponding Sobolev embeddings is lost, and ``critical terms'' are essentially non-compact perturbations of the principal dynamics.  Thus, in general, one would not expect the flow to converge to a compact set, especially in hyperbolic systems. For an overview of associated challenges see the very comprehensive treatise \cite{chu-las-memoirs}. The recent advances in this direction, however, predominantly rely on the dissipation mechanism being supported on the \emph{entire domain}.  
In particular, for structural acoustic interactions, it was shown in \cite{bu-chu-las} that in the presence of an acoustic dissipation distributed over the whole domain domain $\Omega$, is sufficient for existence of a global finite-dimensional attractor, even in presence of a critical acoustic source.

However, the physically relevant localised interior (or boundary) damping poses yet further difficulties, since besides the lack of compactness from criticality, the dissipation now must be ``propagated'' across the domain in order to have  any kind of a stabilizing effect.  A critical source and geometrically restricted damping in a single wave equation had been a long-standing open problem, whose solution ultimately necessitated Carleman estimates and geometric optics analysis \cite{chu-las-tou, chu-las-tou2}. 

In a composite setting the aforementioned techniques would have to  account not merely for two different types of dynamics, but also for the Neumann-type coupling on the interface.  Moreover, the source on the plate  is described by higher-order operators and, if approached via the same strategy as the acoustic counterpart, presents a much more formidable obstacle.  To overcome this difficulty, we develop a new  ``hybrid'' method that employs different techniques applicable to waves with localised damping, while taking advantage of the  full interior dissipation on the plate;  at the same time the analysis accommodates the mixed terms that arise from the coupling.

Ultimately it follows that acoustic damping active \emph{only in a neighbourhood} of the flexible boundary---combined with a nonlinear interior plate damping---still guarantees the existence of a finite-dimensional global attractor, under ``minimal'' assumptions on the nonlinear functions; ``minimal'' here indicating those consistent with the hypotheses which yield a global attractor for the dynamics generated by a single wave equation.

\subsection{Outline of the argument}
To achieve our goal we will appeal to several general criteria from the theory of  infinite-dimensional dynamical systems.
While well-posedness follows via the classical theory of monotone operators in
\cite{barbu}, along with some recent additions devised in \cite{chu-ell-las},
the asymptotic analysis will benefit from novel abstract results gathered in 
\cite{chu-las-memoirs}, specifically tailored for dissipative evolution equations  of {\em hyperbolic} type.
(For general references on dissipative infinite-dimensional systems, the reader is referred to the cornerstones: \cite{babin}, \cite{temam-97} and \cite{hale}; see also \cite{olga} and \cite{chepyzhov-vishik}, the latter addressing non-autonomous systems.)
More specifically, we shall invoke 
\begin{enumerate}[(i)]
\item  a compactness criterion---introduced in \cite{khan} and subsequently recast into a more general abstract form in \cite{chu-las:jde:2007}---recorded here as Theorem~\ref{t:khan} (this result is central to the proof of existence of a global attractor),
\item  a generalization of Ladyzhenskaia's Theorem on dimension of invariant sets, that is Theorem~\ref{t:ladyzhenskaia-general}, which will finally enable us to infer finite-dimensionality of the attractor.
\end{enumerate}

The application of the aforesaid criteria is not immediate, as one might guess in light of the analysis already carried out on the composite PDE system in \cite{bu-chu-las} with (acoustic) {\em full interior} damping, and of the novel tools devised in \cite{chu-las-tou} for  the (uncoupled) wave equation with {\em localised} dissipation.
A primary source of difficulty stems from the fact that these abstract results require  estimates on the quadratic energy corresponding to the {\em difference of two trajectories}, rather than to a {\em single trajectory}, and the presence of localised acoustic dissipation in conjunction with the coupling
(with the plate equation) brings about further technical difficulties, over those encountered in \cite{bu-chu-las} and \cite{chu-las-tou}, as explained in 
\S\ref{ss:newchallenges}.

The analysis to be carried out includes various steps. After a preliminary discussion of well-posedness of the PDE system in an appropriate  state space, a major task will be to establish {\em asymptotic compactness}.  That the global attractor has a finite fractal dimension will be obtained subsequently, through a use of both a Carleman-like identity established in \cite{chu-las-tou} for the wave dynamics, and exploiting the attractor's compactness. As  expected, proving finite dimensionality and the additional smoothness of the attractor will require stronger assumptions on the nonlinear dissipation terms in either of the coupled equations.

Note that we chose to separate the proof of existence of a global attractor from that of its finite-dimensionality, as the former property does hold under weaker assumptions. In this respect, the present argument differs from the ones followed in a large part of the recent literature on dissipative dynamical systems, where existence of an  {\em exponential}\footnote{The concept of exponential attractor was introduced in \cite{eden-foias}.} attractor is sought, as it brings about---besides a certain robustness under perturbations and numerical approximations---the {\em intrinsic} feature of finite dimensionality; see, e.g. (and without any claim of completeness), \cite{eden-milani}, \cite{efendiev-miranville}, \cite{fabrie-etal}, \cite{grasselli-pata}, \cite{grasselli-prazak}, 
\cite{schimperna}.
The question whether the PDE system under investigation possesses an exponential attractor will not be discussed here.

We finally note that the boundedness of the attractor in a smoother functional space is closely tied to its finite-dimensionality. The proofs of both results are interconnected; it is precisely the compactness of the attractor in the finite-energy space and the regularity estimate in the \emph{higher-energy} space which together show that critical (non-compact) sources, do not substantially perturb the structure of the attracting set. A recent result \cite{conti-pata} provides an elegant way to prove regularity of a global attractor in higher norms, without directly appealing to uniform-in-time estimates; however, the application of this argument to a wave equation requires a stronger damping than employed in the present case.

\smallskip

A brief outline of the paper follows.
\begin{itemize}
	\item  
In Section~\ref{s:pde-model} we introduce the PDE model under investigation, 
along with the statements of our main results: well-posedness (Theorem~\ref{thm:wellposed}), existence of a global attractor for the corresponding dynamics (Theorem~\ref{thm:attractor-existence}), and finite-dimensionality, as well as  regularity of the attractor (Theorem~\ref{thm:finite-dim}).  
The natural energies and state spaces for either uncoupled equation, and then
for the coupled PDE system, are introduced and briefly discussed in Section~\ref{subsec:energies}.

\item Section~\ref{s:introductory-results} contains some preliminary work for the proof of existence of the attractor. More precisely, we establish a basic identity involving the
integral of the quadratic wave and plate energies pertaining to 
the {\em differences of trajectories}; see Lemma~\ref{l:full-basic-identity}.

\item Section~\ref{s:existence-attractor} is devoted to the proof of Theorem~\ref{thm:attractor-existence}.
Existence of a global attractor is established in view of the gradient structure of the dynamical system and {\em asymptotic smoothness} of the semi-flow.  In turn, this latter (crucial) property follows from a {\em pointwise} estimate of the total quadratic energy (Proposition~\ref{p:energy-pointwise-estimate}), combined with a suitable ``weak sequential compactness'' property (Proposition~\ref{p:compact}).  
\item Section~\ref{s:finite-dim} is focused on the proof of Theorem~\ref{thm:finite-dim},
technically based on the key inequality established in Lemma~\ref{l:central-to-finitedim}
(refined in Lemma~\ref{p:pre-final}), and on the Proposition~\ref{p:forces}. Finally, a short Appendix concludes the paper.
\end{itemize}

\section{The PDE model and the statement of main results} \label{s:pde-model}
The PDE system under investigation models structure-acoustics interactions described by the acoustic velocity potential $\zet(x,t)$ evolving in a smooth domain 
$\Om\subset \bbR^3$, and the vertical displacement $w(x,t)$ of a thin hinged plate whose midsection occupies a part of the boundary of $\Om$; we denote the latter by 
$\Gam:=\p\Om$, and assert that it consists of two relatively open simply-connected segments with overlapping closures $\Gam_0\cup \Gam_1$, where $\Gam_0$ represents the plate, and $\Gam_1$ is the rigid wall of the acoustic chamber. The acoustic damping is confined to a \emph{subset} of $\Om$ corresponding to the support of the cutoff function $\chi(x)$.

\begin{equation}\label{e:pde-system-0}
\left\{
\begin{array}{lll}
\zeta_{tt}-\Delta \zeta + \chi(x) g(\zeta_t) = f(\zeta)
& {\rm in} & \Omega\times (0,T)
\\[2mm]
\displaystyle \frac{\partial \zeta}{\partial\nu}  = 0
& {\rm on} & \Gamma_1\times (0,T)
\\[2mm]
\displaystyle\frac{\partial \zeta}{\partial\nu} = \alpha \kappa \,w_t
& {\rm on} & \Gamma_0\times(0,T)
\\[2mm]
\zeta(0,\cdot) = \zeta_0\,, \quad \zeta_t(0,\cdot) = \zeta_1 
& {\rm in} & \Omega
\\[2mm]
w_{tt}+\Delta^2 w +  b(w_t) + \beta \kappa \zeta_t|_{\Gamma_0} = f_2(w)
& {\rm in} & \Gamma_0\times(0,T)
\\[2mm]
w=\Delta w = 0 & {\rm on} & \partial\Gamma_0\times(0,T)
\\[2mm]
w(\cdot,0) = w_0\,, \quad w_t(\cdot,0) = w_1 & {\rm in} & \Gamma_0\,.
\end{array}
\right.
\end{equation}
The positive constants $\al$, $\bet$ and $\kap$ represent parameters dictated by the physical model in question. Following \cite{bu-chu-las}, we shall consider, specifically, 
\begin{equation}\label{e:berger}
f_2(w) = p_0 - \Big[Q-\int_{\Gamma_0} |\nabla w|^2\Big]\Delta w\,,
\end{equation}
which is the semilinear term that occurs in {\em Berger's plate model}:
the function $p_0$, which is related to transversal forces, belongs
to $L_2(\Gamma_0)$, while the real constant $Q$ describes in-plane  forces applied to the plate; for more on modeling according to the Berger approach, see \cite{chueshov-book}.

\begin{remark}
\begin{rm}
It is important to emphasise that although we consider a plate equation with
a nonlinear term $f_2$ of the form \eqref{e:berger}, the analysis carried 
out in the present paper might be extended to more general nonlinear functions, 
provided they satisfy proper conditions, such as those listed in Assumption~4.11 
in \cite{chu-las-memoirs} (following this reference's notation, in the present
case $\sA$ would denote the plate dynamics operator, that is the realization of the bilaplacian $\Delta^2$ with hinged boundary conditions).
\end{rm}
\end{remark}

\subsection{Assumptions and main results}
\begin{figure}[h]
\includegraphics[width=2.4in,keepaspectratio=true]{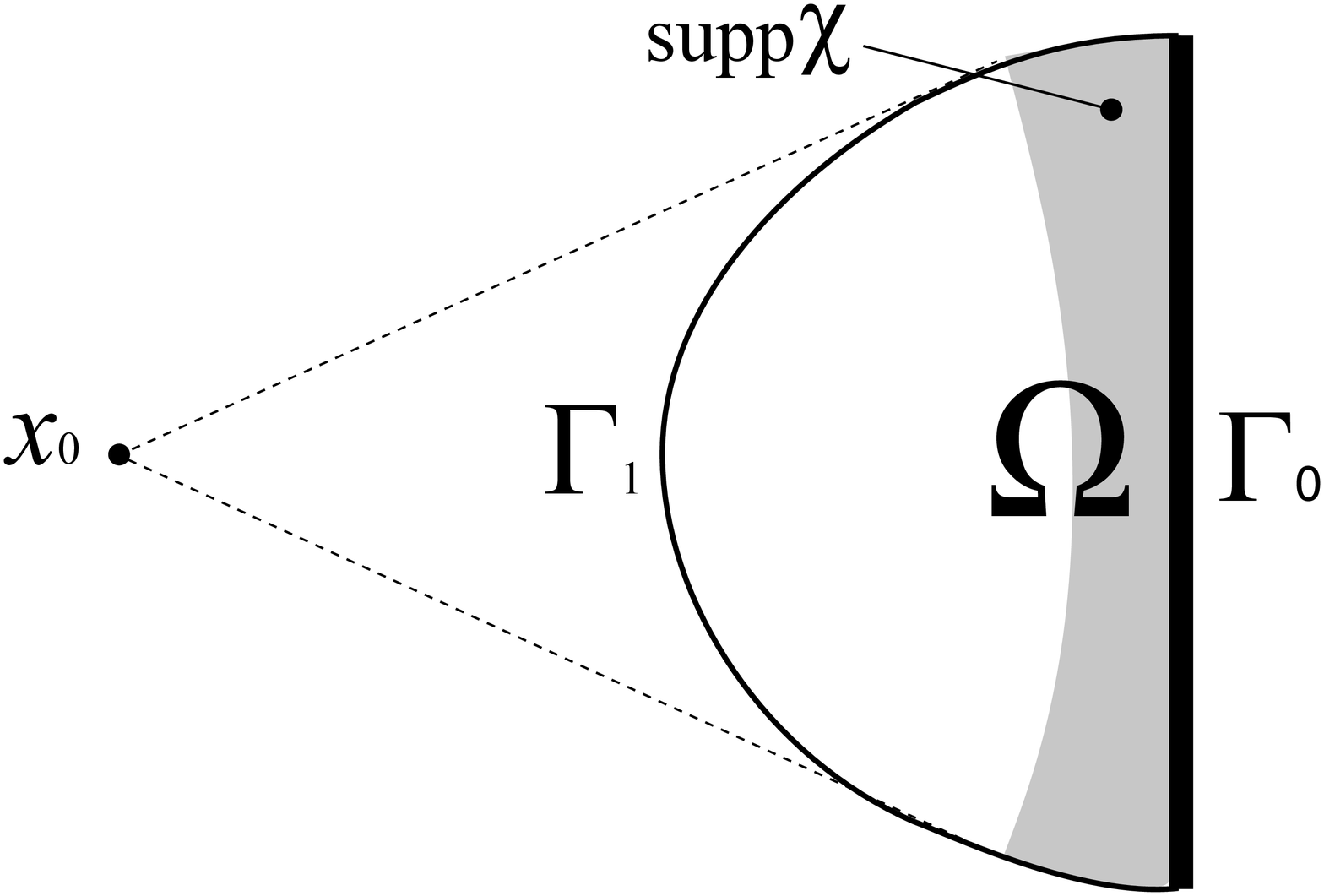}
\caption{\label{fig:domain}
A cross-section of a suitable domain $\Omega$. The rigid wall $\Gam_1$, and, in particular, the part of the boundary outside $\mbox{supp}\chi$ is convex; however, more general configurations are possible: for other examples see \cite{las-trig-zhang}.
}
\end{figure}

Due to geometric restrictions on the acoustic damping, the shape of $\Om$ must satisfy certain restrictions in order for the feedback $\chi(x) g(\zet_t)$ to be ``effective,'' as dictated by geometric optics and the classical results on propagation of singularities   \cite{bar-leb-rauch}.  A sufficient assumption is the following:
\begin{assumption}[Geometry of the domain]
\label{a:geometry}
$\Gamma_1$ is a level-surface of a convex function and there exists $x_0\in \mathbb{R}^n$ such that 
\begin{equation*}
(x-x_0)\cdot \nu \le 0 \quad \textrm{on } \; \Gamma_1\,,
\end{equation*}
with $\nu$ being the outward normal vector field on $\Gamma_1$.
\end{assumption}
Figure \ref{fig:domain} shows what a cross-section of a possible domain may look like.

\begin{remark}[More general geometric conditions]
\begin{rm}
The convexity of $\Gam_1$ is a sufficient requirement, but it is possible to further 
relax this assumption. In particular, it suffices to ask for the subset of $\Gam_1$ 
that lies away from $\operatorname{supp} \chi(x)$ to be a level-surface of a function 
with a non-vanishing gradient and positive-definite or negative-definite Hessian (in 
the latter case the condition on $x_0$ should be changed to $(x-x_0) \cdot \nu \ge 0$). 
See \cite{las-trig-zhang} for more details and examples.
\end{rm}
\end{remark}

\begin{assumption}[Nonlinear terms]  \label{a:sources}
The function $f$ satisfies:
\begin{itemize}
\item
$f\in C^2(\mathbb{R})$ and there exists $C_f$ such that
\begin{equation}\label{h:source-1}
|f''(s)|\le C_f (1+|s|) \qquad \forall s\in\bbR;
\end{equation}
\item
$\displaystyle\limsup_{|s|\to \infty}\frac{f(s)}{s}=:-2\lambda <0$.
\end{itemize}
Accordingly, set $f_1(s) = f(s) + \lambda s$.
\end{assumption}

Existence of a global attractor will be established under rather weak
assumptions on the feedback maps $g$ and $b$.

\begin{assumption}[Damping terms, I] \label{a:dampings-0}\leavevmode
\begin{enumerate}
\item
$g\in C(\mathbb{R})$ is a monotone increasing function such that $g(0)=0$; 
in addition, $g$ is linearly bounded at infinity, i.e. there exist positive
constants $m_g$, $M_g$ such that 
\begin{equation}\label{h:wdamping}
m_g|s|\le |g(s)| \le M_g|s| \qquad \textrm{for \;} |s|\ge 1\,.
\end{equation}
\item
$b\in C(\mathbb{R})$ is a monotone increasing function such that 
$b(0)=0$; and there exists $m_b$ such that 
\begin{equation}\label{h:pdamping}
|b(s)|\ge m_{b} |s| \qquad \textrm{for \;} |s|\ge 1\,.
\end{equation}

\end{enumerate}
\end{assumption}

\smallskip
\noindent
A few comments about the above assumptions are in order.
As observed already in \cite{bu-chu-las} and \cite{chu-las-tou}, a prime consequence of Assumption~\ref{a:sources} is local Lipschitz continuity of  the nonlinear term $f$ (and of $f_1$, as well), as an operator from $H^1(\Omega)$ into  $L_2(\Omega)$. Namely, one has 
\begin{equation*}
\|f(\zeta_1)-f(\zeta_2)\|\le C(r) \|\zeta_1-\zeta_2\|_{1,\Omega} \qquad \forall \zeta_i 
\quad \textrm{with}\; \|\zeta_i\|_{1,\Omega}\le r\,, \, i=1,2\,.
\end{equation*}
The above property is easily shown by using well known Sobolev embedding results
in $3$-D domains, which establish the {\em critical} threshold in the polynomial 
bound on $f$.
The same (Lipschitz continuity) is true for the nonlinear term $f_2$ which occurrs in the plate equation; more precisely,
\begin{equation*}
\|f_2(w_1)-f_2(w_2)\|\le C(r) \|w_1-w_2\|_{2,\Gamma_0} \qquad \forall w_i 
\quad \textrm{with}\; \|w_i\|_{2,\Gamma_0}\le r\,, \, i=1,2\,.
\end{equation*} 

Next, observe that the basic Assumption~\ref{a:dampings-0} forces the acoustic damping to be linearly bounded at infinity. It is important to emphasise that this condition is actually necessary, because of 
the geometrical restrictions on the damping term. This fact was first exhibited in \cite{vancoste-martinez} for a wave equation with boundary dissipation; indeed, the case of localised dissipation around a portion of the boundary yields the same technical difficulties and the same results as that of boundary damping.
The reader is referred to \cite[Section~1]{chu-las-tou} for a more detailed discussion of these issues.

The aforementioned basic assumptions are sufficient to establish both\break well-posedness for 
the PDE problem \eqref{e:pde-system-0} and existence of a global attractor for the corresponding dynamical system.

\begin{theorem}[Well-posedness]\label{thm:wellposed}
Define the state space 
\[
\sY\dfn H^1(\Om) \times L^2(\Om) \times H^2_0(\Gam_0)\times L^2(\Gam_0)\,.
\]
The system \eqref{e:pde-system-0} generates a strongly continuous semigroup $S(t)$ on $\sY$. In particular, given initial data, at $t=0$, $W_0=\{\zet_0, \zet_1, w_0, w_1\}\in \sY$ there exists a unique solution to \eqref{e:pde-system-0}
\[
S(\cdot)W = \{\zet,\zet_t,w,w_t\}\in C([0,T],\sY)
\]
for any $T\geq 0$. The solution  satisfies the following variational identities
\begin{equation}\label{t:wellpos:e:wave}
\begin{split}
(\zet_t, \phi)\bigg|_0^T - \intT (\zet_t, \phi_t)  + \intT (\grad \zet, \grad \phi) 
+& \intT (\chi g(\zet_t),\phi) \\
=& \intT (f(\zet),\phi) +\intT \langle \al\kap w_t, \phi\rangle\,,
\end{split}
\end{equation}

\begin{equation}\label{t:wellpos:e:plate}
\begin{split}
\lg w_t, \psi\rg\bigg|_0^T - \intT \lg w_t, \psi_t\rg  + \intT \lg \Del w, \Del \psi\rg 
 +& \intT \lg  b(w_t),\psi \rg \\
=& \intT \lg f_2(w),\psi\rg -\intT \lg \psi, \bet\kap \zet_t \rg\,,
\end{split}
\end{equation}
for any test functions
\[
\phi\in C^1([0,T],L^2(\Om))\cap C([0,T]; H^1(\Om))\,,\;
\psi\in C^1([0,T],L^2(\Gam_0))\cap C([0,T]; H^2_0(\Gam_0))\,.
\]
Moreover, if the initial data belong to the spaces
\[
\{\zet_0, \zet_1\} \in H^2(\Om)\times H^1(\Om),
\quad \{w_0,w_1\}\in\big(H^4(\Gam_0)\cap H_2(\Gam_0)\big)\times H^2_0(\Gam_0)
\]
with the compatibility conditions
\[
\Dn{\zet_0}\Big|_{\Gam_1} = 0,\qquad  \Dn{\zet_0}\Big|_{\Gam_0} = \al\kap w_1\,,
\]
then
\[
\{\zet,\zet_t,w,w_t\} \in L^\infty\big(0,T;H^2(\Om)\times H^1(\Om)
\times H^4(\Gam_0)\times H^2(\Gam_0)\big)\,.
\]
\end{theorem}

\medskip
We introduce the set $\sN$ of equilibria for the dynamical system $(\sY,S(t))$;
namely
\begin{equation} \label{e:equilibria}
\sN := \{W \in \mathcal{Y}\;:\; S(t)W =  W\quad \text{for all} \quad t\geq 0\}\,.
\end{equation}

\begin{theorem}[Existence of the global attractor] \label{thm:attractor-existence}
Suppose that the Assumptions  \ref{a:geometry}--\ref{a:sources} hold.
If the damping functions satisfy the Assumption~\ref{a:dampings-0},
then the dynamical system $(\sY,S(t))$ generated by the PDE problem
\eqref{e:pde-system-0} has a compact global attractor $\mathscr{A}\subset \sY$, 
which coincides with the unstable manifold $M^u(\sN)$ of the set $\sN$ of 
stationary points:
\begin{itemize}
\item 
$\mathscr{A}\equiv M^u(\sN)$;
\item 
$\displaystyle\lim_{t\to +\infty}{\rm dist}(S(t)W,\sN)=0\quad
\textrm{for any}\quad W \in \sY$.
\end{itemize}
\end{theorem}

To prove that the attractor has a finite fractal dimension we will
need to strengthen the regularity and growth condition on  both damping functions.

\begin{assumption}[Dampings, II] \label{a:dampings-1}
Assume $g,b\in C^1(\mathbb{R})$, with $g(0)=0=b(0)$, and 
\begin{enumerate}
\item
there exist positive constants $m_g'$, $M_g'$ such that 
\begin{equation}\label{e:wave-hypo-stronger}
m_g'\le |g'(s)| \le M_g' \qquad \forall s\in \mathbb{R}\,,
\end{equation}
\item
there exist positive constants $m_b'$, $M_b'$ such that 
\begin{equation}\label{e:plate-hypo-stronger}
m_b'\le |b'(s)| \le M_b'(1+sb(s)) \qquad \forall s\in \mathbb{R}\,.
\end{equation}
\end{enumerate}

\end{assumption}
Then, the following result holds.

\begin{theorem}[Finite dimensionality and regularity of the attractor] 
\label{thm:finite-dim}
Suppose the hypotheses of Theorem~\ref{thm:attractor-existence} are satisfied.
If, in addition, the Assumption~\ref{a:dampings-1} holds, 
then the global attractor $\mathscr{A}$ has a finite fractal dimension.

Furthermore, the attractor $\mathscr{A}$ is bounded in the domain of the nonlinear semigroup generator; in particular, $\mathscr{A}$ is a bounded subset of 
\begin{equation*}
H^2(\Omega)\times H^1(\Omega)\times H^4(\Gamma_0)\times H^2(\Gamma_0)\,.
\end{equation*}
\end{theorem}


\subsection{Energy identity and bounds} \label{subsec:energies}
With the system \eqref{e:pde-system-0} we associate the following energy functionals:
\begin{eqnarray*}
\sE_\zet(t) &\dfn& \frac12 \left( \|\grad \zet\|^2 + \lambda\|\zet\|^2+\|\zet_t\|^2\right)- \intOm F_1(\zet(t))\\
\sE_w(t) &\dfn& \frac12\Big(\|\Del w\|^2 +  \frac{1}{2}\|\grad w\|^4+ \|w_t\|^2 \Big) 
- \Pi(w(t))\,,
\end{eqnarray*}
where $F_1(s)$ is the antiderivative of $f_1(s)$ vanishing at $0$, and 
\begin{equation} \label{e:plate-functional}
\Pi(w) =  (p_0, w) + \frac{Q}{2} \|\grad w\|^2  - \frac{1}{4}\|\grad w\|^4\,.
\end{equation}
Define  the total energy
\begin{equation}\label{d:total-energy}
\sE(t) = \sE_{\zet,w}(t)\dfn \beta \sE_\zet(t) + \alpha\sE_w(t)\,.
\end{equation}
The following identity satisfied by $\sE(t)$ is standard due to the fact that the  dissipation feedbacks are linearly bounded at infinity and the sources correspond to  locally Lipschitz operators on the energy space. The equation can be derived for strong solutions by using test functions $\bet\zet_t$ and $\al w_t$ in the variational identities \eqref{t:wellpos:e:wave} and \eqref{t:wellpos:e:plate} respectively. 
Since the result is continuous with respect to the finite-energy topology, it can be extended to all weak solutions: 
\begin{equation} \label{e:energy-identity}
\sE(t) + \beta \int_s^t (\chi g(\zeta_t),\zeta_t) + \alpha
\int_s^t(b(w_t),w_t)= \sE(s)\,.
\end{equation}
Also, let us introduce positive quadratic energy functionals:
\[
E_\zet(t) \dfn \frac12\big( \|\grad \zet\|^2 + \lambda\|\zet\|^2 + \|\zet_t\|^2\big)\,,\quad
E_w(t) \dfn \frac12\Big(\|\Del w\|^2+\frac{1}{2}\|\grad w\|^4+\|w_t\|^2 \Big)\,;
\]
\begin{equation}\label{d:pos-energy}
E(t) = E_{\zet,w}(t) \dfn E_\zet(t) + E_w(t)\,.
\end{equation}

Owing to the dissipativity property (in Assumption~\ref{a:sources}) satisfied by $f$ and to the structure of the functional $\Pi$ in \eqref{e:plate-functional}, it is not difficult to obtain upper and lower bounds for the full energy of the system;
see also \cite[Section~2.2]{bu-chu-las}.
Explicit computations pertaining to the wave component are found in 
\cite[Section~2]{chu-ell-las}.

\begin{proposition}[Bounds on the energy]\label{p:energy-bound}
Let $\sB$ be a bounded subset of $\sY$. If $\{\zet_0,\zet_1,w_0,w_1\} \in \sB$ then there exists constants $C_{1,\sB}, C_{2,\sB}$ dependent only on the diameter of $\sB$ (in the topology of  $\sY$) such that
\begin{equation}\label{e:bounds}
C_{1,\sB} E(t) - C_{2,\sB} \leq  \sE(t) \leq \sE(0) \qquad \forall t\geq 0\,.
\end{equation}
\end{proposition}

\medskip
We conclude this section by introducing the abstract dynamic operators pertaining
to either equation, namely:
\begin{align*}
& A_1  v = -\Delta v,\quad 
\sD(A_1) = \Big\{  v \in H^2(\Om)\,:\;   \Dn{v} \Big|_{\Gam} =0\Big\}\,;
\\
& A_2  v = \Delta^2 v,\quad \sD(A_2) = \Big\{ v \in H^4(\Gam_0)\cap H^1_0(\Gam_0)\,:
;  \Del v \big|_{\p\Gam_0} =0\Big\}\,.
\end{align*}


\section{The differences of trajectories: introductory results}
\label{s:introductory-results}
Seeking to apply the abstract results recorded in the Appendix as Theorem~\ref{t:khan} and Theorem~\ref{t:ladyzhenskaia-general}, to investigate the asymptotic behaviour of the solutions to \eqref{e:pde-system-0}, we must study  {\em differences of its trajectories}, rather than the trajectories themselves.
In this section we introduce the relative basic definitions, along with a series 
of preliminary identities which constitute a first step in the proof of our main results.

\subsection{Auxiliary functions and parameters}
Given two different evolution trajectories $(h,h_t,u,u_t)$ and 
$(\zeta,\zeta_t,w,w_t)$ of the coupled PDE system \eqref{e:pde-system-0},
we introduce the differences $z$ and $v$, namely,
\begin{align}
z&:= h-\zeta\,,  \qquad(\zeta =h+z\,); 
\\
v&:= u-w\,,  \qquad(u =w+v\,).
\end{align}
The pair $(z,v)$ readily solves the new coupled system 
\begin{equation}\label{e:pde-system-1}
\left\{
\begin{array}{lll}
z_{tt}-\Delta z +\lambda z + \chi(x) \tg(\zeta_t) = \tfw(z)
& {\rm in} & \Omega\times (0,T)
\\[2mm]
\displaystyle \frac{\partial z}{\partial\nu}  = 0
& {\rm on} & \Gamma_1\times (0,T)
\\[2mm]
\displaystyle\frac{\partial z}{\partial\nu} = \alpha \kappa \,v_t
& {\rm on} & \Gamma_0\times(0,T)
\\[2mm]
z(0,\cdot) = h_0-\zeta_0\,, \quad z_t(0,\cdot) = h_1-\zeta_1 
& {\rm in} & \Omega
\\[2mm]
v_{tt}+\Delta^2 v +  \tb(w_t) + \beta \kappa z_t|_{\Gamma_0} = \tfp(v)
& {\rm in} & \Gamma_0\times(0,T)
\\[2mm]
v=\Delta v = 0 & {\rm on} & \partial\Gamma_0\times(0,T)
\\[2mm]
v(0,\cdot) = u_0-w_0\,, \quad v_t(0,\cdot) = u_1-w_1 & {\rm in} & \Gamma_0\,.
\end{array}
\right.
\end{equation}
where we have set
\begin{align}
\tg(z_t)&:= g(\zeta_t+z_t)-g(\zeta_t)\,, 
\label{d:gtilde}\\
\tb(v_t)&:= b(w_t+v_t)-b(w_t)\,,
\label{d:btilde}\\
\tfw(z)&:= f_1(\zeta+z)-f_1(\zeta)\,, 
\label{d:f1tilde}\\
\tfp(v)&:= f_2(w+v)-f_2(w)\,.
\label{d:f2tilde}
\end{align}
Technically each of the introduced functions also depends on one of the terms in the corresponding difference, but that fact will be suppressed for brevity of notation.

\subsection{Smooth cutoff functions}\label{sec:cutoffs}
Following \cite[Section~6]{toundykov-na07}, we introduce two smooth cutoff 
functions whose role is to single out in $\Omega$ the dissipative and the non-dissipative subdomains.
More precisely, $\psi, \phi: \overline{\Omega}\to [0,1]$ are $C^2$ functions 
with the following properties:
\begin{itemize}
\item
${\rm supp}(\phi)\subset {\rm supp}(\chi)$,
\item
$\Omega \setminus {\rm supp}(\chi)\subset \{x:\, \psi(x)=1\}$,
\item
$\psi \equiv 0$ in a neighbourhood of $\Gamma_0$, 
\item
for any $x\in \overline{\Omega}$ at least one of $\psi(x)$ and $\phi(x)$ 
equals $1$, namely 
\begin{equation*}
\{ x\in \overline{\Omega}: \psi(x)<1\}\subset \phi^{-1}(\{1\})\,.
\end{equation*}
\end{itemize}
By setting $\hz:= \psi z$ and $\cz:= \phi z$, it is easily verified that $\hz$ and $\cz$ satisfy 
\begin{equation}
\hz_{tt} -\Delta \hz +\lambda\hz + [\![\Delta,M_\psi]\!]z + \psi \chi \tg(z_t)
=\psi\tfw(z)\,,
\end{equation}
\begin{equation}
\cz_{tt} -\Delta \cz +\lambda\cz + [\![\Delta,M_\phi]\!]z + \phi \chi \tg(z_t)
=\phi\tfw(z)\,,
\end{equation}
where $M_\psi$, $M_\phi$ denote pointwise (a.e.) multiplication by $\psi$ or 
$\phi$, respectively, while the pairing $[\![\cdot,\cdot]\!]$ denotes a commutator,
acting as follows:
\begin{equation*} 
[\![\Delta,M_\psi]\!]z:= \Delta(\psi z)- \psi\Delta z\,.
\end{equation*}
A straightforward computation gives the equivalent form 
\begin{equation}\label{e:commutator}
[\![\Delta,M_\psi]\!]z= \Delta\psi\, z+2 \nabla\psi\nabla z\,.
\end{equation}

\subsection{Vector field $\vh$}\label{s:vec-field}
Let us recall from \cite{chu-las-tou} the construction of a function $d$ and 
the vector field 
\begin{equation} \label{h-vector}
\vh(x) :=\nabla d(x)\,,
\end{equation}
whose key properties are 
\begin{enumerate}[(i)]
\item
\begin{equation}
\min_{x\in \overline{\Omega}}d(x)>0\,;
\end{equation}
\item
\begin{equation}
(\vh\cdot\nu)|_{\Gamma_1}=\frac{\partial d}{\partial \nu}\Big|_{\Gamma_1} \equiv 0\,;
\end{equation}

\item
the Jacobian matrix $J_{\vh}$ of $\vh$---which coincides with the 
hessian matrix $H_d$ of $d$---evaluated on $\Gamma_1$ is positive
definite. In particular, $d$ can be extended to some open set 
containing all of $\bar{\Omega}$ so that for some $\rho>0$ one has
\begin{equation} \label{e:positive-hessian}
H_d= J_{\vh} \ge \rho\,I
 \quad \forall x\in \Omega\,.
\end{equation}
The reader is referred to \cite[p.~301-303]{las-trig-zhang} for more details.
\end{enumerate}
The (Carleman-type) estimates pertaining to the wave component of the system
will involve the pseudo-convex function 
$\Phi: \Omega\times \mathbb{R}\to \mathbb{R}$ defined by
\begin{equation} \label{d:pseudo-convex}
\Phi(x,t) := d(x) - c\big|t-T/2\big|^2\,,
\end{equation}
where at the outset $T>0$ and $c$ is a non-negative constant, with 
$T$ large enough to satisfy
\begin{equation}\label{e:constraint-2}
T> 2\big(\displaystyle\max_{x\in \overline{\Omega}}d(x)/c\big)^{1/2}\,,
\end{equation}
thus ensuring $\Phi(x,0), \Phi(x,T)<0$.

\subsection{Preliminary fundamental identities}

\subsubsection{Wave component}
Our starting point is a key identity pertaining to the wave component
of the PDE system \eqref{e:pde-system-1} satisfied by the differences 
$(z,v)$. 
This result has been established in \cite{chu-las-tou}; see
Proposition~5 in \S6.4 therein.
Because of the slightly different wave energy, in the present case 
the identity reads as follows.


\begin{proposition}[Wave Fundamental identity (Carleman-type), \cite{chu-las-tou}]
\label{p:fundamental-identity}
Suppose that the Assumptions~\ref{a:geometry}, \ref{a:sources}, \ref{a:dampings-0} hold.
Take smooth initial data  $(\zeta_0,\zeta_1), (h_0,h_1)\in \sD(A_1)$, and set $z=h-\zeta$.
Let $\Phi(x,t)$ be given by \eqref{d:pseudo-convex}, and let 
$\vh :=\nabla \Phi = \nabla d$. 
Recall the notation $\hz:= \psi z$ and $\cz:= \phi z$
($\psi$, $\phi$ are the cutoff functions introduced in Section \ref{sec:cutoffs}). 
Then, for any $\tau \ge 0$ and any positive constant $\bsC$ one has 

\begin{equation}\label{e:fundamental-carleman}
\addtolength{\fboxsep}{5pt}
\boxed{
\begin{split} 
\int_{Q_T} e^{\tau\Phi} & (J_{\vh}-\rho Id)\nabla \hz \cdot \nabla \hz 
+ [Energy]_\psi + \bsC\, [Energy]_\phi
\\
& = [Damping]- \tau \int_{Q_T}e^{-\tau\Phi}\sM_1^2
+ \int_{Q_T}\psi\big(\tfw(z)-\chi\tg(z_t)\big)\sM_1
\\
& \quad -\int_{Q_T}([\![\Delta,M_\psi]\!]z)\,\sM_1 + [\text{Almost lower order}]
\\
& \quad +\,(BT)_{\Sigma} - \widehat{C}_{0,T}
- \bsC\,\overset{\circ}{C}_{0,T}\,,
\end{split}
}
\end{equation}
with $\sM_1$ and $\mu$ defined by 
\begin{subequations}\label{e:multipliers}
\begin{eqnarray}
& & \sM_1 := e^{\tau\Phi}\big(\vh\cdot \nabla \hz -\Phi_t\hz_t\big)\,,
\\[1mm]
& & \mu := {\rm div}\,(e^{\tau\Phi}\vh)- \partial_t(e^{\tau\Phi}\Phi_t)\,,
\end{eqnarray}
\end{subequations}
and 
\begin{eqnarray}
& & [Energy]_\psi := \left(\frac{\rho}2-c\right)
\int_{Q_T}e^{\tau\Phi}\big(|\nabla\hz|^2+\lambda\hz^2+{\hz}_t^2\big)\,,
\\
& & [Energy]_\phi := 
\int_{Q_T}e^{\tau\Phi}\left(|\nabla\cz|^2+\lambda\cz^2+{\cz}_t^2\right)\,;
\end{eqnarray}
moreover, we set 
\begin{equation}
\begin{split}
[Damping] :=\, &
2\,\bsC \int_{Q_T} e^{\tau\Phi}\cz_t^2 - \int_{Q_T}\psi\chi \tg(z_t)\hz
\left[\frac{\mu}2-\big(\frac{\rho}2+c\big)e^{\tau\Phi}\right]
\\
& \quad -\,\bsC\int_{Q_T}\phi\chi\tg(z_t)e^{\tau\Phi}\cz\,,
\end{split}
\end{equation}
and
\begin{eqnarray}
\lefteqn{
[\text{Almost lower order}] := 
- \int_{Q_T} \hz\nabla\hz\cdot \nabla
\Big(\frac{\mu}2-\big(\frac{\rho}2+c\big)e^{\tau\Phi}\Big)}
\notag \\
& & 
+ \int_{Q_T} \hz\hz_t\frac{\partial}{\partial t}
\Big(\frac{\mu}2-\big(\frac{\rho}2+c\big)e^{\tau\Phi}\Big)
+ \bsC\int_{Q_T} \cz\cz_t\frac{\partial}{\partial t}e^{\tau\Phi}
- \bsC \int_{Q_T} \cz\nabla\cz\cdot \nabla e^{\tau\Phi}
\notag \\
& & -\int_{Q_T}\Big([\![\Delta,M_\psi]\!]z-\psi\tfw(z)\Big)
\Big\{\sM_1+ \hz\Big[\frac{\mu}2-\big(\frac{\rho}2+c\big)e^{\tau\Phi}\Big]
\Big\}
\notag \\
& & 
-\lambda\int_{Q_T}\hz \Big[\sM_1+ \hz\big(\frac{\mu}2+\rho e^{\tau\Phi}\big)
\Big]
\notag \\
& & 
-\bsC\int_{Q_T}\big(\,[\![\Delta,M_\phi]\!]z-\phi\tfw(z)\,\big)
e^{\tau\Phi}\cz\,.
\end{eqnarray}
The boundary terms are collected in $(BT)_{\Sigma_T}$:
\begin{equation}
\begin{split}
(BT)_{\Sigma_T} := \,& 
\int_{\Sigma_T} \frac{\partial \hz}{\partial\nu} \sM_1
+ \int_{\Sigma_T} \frac{\partial \hz}{\partial\nu} \hz 
\left[\frac{\mu}2-\big(\frac{\rho}2+c\big)e^{\tau\Phi}\right]
\\
&  +\,\frac12\,\int_{\Sigma_T} e^{\tau\Phi}
\big(\hz_t^2-|\nabla\hz|^2\big)(\vh\cdot\nu)\,.
+ \bsC \int_{\Sigma_T} \frac{\partial \cz}{\partial\nu} \cz e^{\tau\Phi}\,;
\end{split}
\end{equation}
while $\widehat{C}_{0,T}$ and $\overset{\circ}{C}_{0,T}$
are given, respectively, by:
\begin{subequations} \label{e:time-traces}
\begin{align} 
\widehat{C}_{0,T} & := 
\Big[\int_{\Omega}e^{\tau\Phi}\hz_t (\vh\cdot \nabla\hz)\Big]_0^T 
-\frac12\Big[\int_{\Omega}e^{\tau\Phi}\Phi_t(\hz_t^2+|\nabla\hz|^2)\Big]_0^T 
\notag \\
& \qquad \qquad + \,\Big[\int_{\Omega}\hz_t\hz
\Big(\frac{\mu}2-\big(\frac{\rho}2+c\big)e^{\tau\Phi}\Big)\Big]_0^T\,,
\\[2mm]
\overset{\circ}{C}_{0,T}&:=
\Big[\int_{\Omega}\cz_t\cz e^{\tau\Phi}\Big]_0^T\,.
\end{align}
\end{subequations}

\end{proposition}

\begin{remark}
\begin{rm}
The proof of Proposition~\ref{p:fundamental-identity} is rather technical; it involves successive application of weighted multipliers:  $\sM_1$, $\sM_2=\mu\hz$ ($\sM_1$ and $\mu$ are defined 
in \eqref{e:multipliers}) and $\sM_3=e^{-\tau\Phi}\cz$.
The reader is referred to \cite[Section~7.3]{chu-las-tou} for all the details.
\end{rm}
\end{remark}

Rewriting the fundamental identity \eqref{e:fundamental-carleman} for
$c, \tau =0$ or, equivalently, substituting the identity $I$ in place of 
$e^{\tau \Phi}$ and utilizing the `simplified' multipliers
\begin{equation*}
\sM_1 = \vh\cdot \nabla \hz\,, \qquad
\mu={\rm div}\,\vh\,,
\end{equation*}
we obtain the following basic assertion, which is
Proposition~3 in \cite[\S5.1]{chu-las-tou}.

\begin{lemma}[Wave Basic Identity]
\label{l:basic-identity}
Suppose that the Assumptions~\ref{a:geometry}, \ref{a:sources}, \ref{a:dampings-0} hold.
Take smooth initial data $(\zeta_0,\zeta_1), (h_0,h_1)\in \sD(A_1)$, 
and set $z=h-\zeta$.
Recall the notation $\hz:= \psi z$ and $\cz:= \phi z$
($\psi$, $\phi$ are the cutoff functions introduce in Section \ref{sec:cutoffs}), and the vector field
$\vh$.
Then, for any positive constant $\bsC$ one has 
\begin{equation}\label{i:wave-basic}
\begin{split} 
\int_{Q_T} (J_{\vh}-\rho I)& \nabla \hz \cdot \nabla \hz 
+ 
\frac{\rho}2\,\int_{Q_T}\big(|\nabla\hz|^2+\lambda\hz^2+{\hz}_t^2\big)
+ \bsC\,\int_{Q_T}\big(|\nabla\cz|^2+\lambda\cz^2+{\cz}_t^2\big)
\\
& = 2\,\bsC \int_{Q_T}\cz_t^2 + \lambda \rho \int_{Q_T}\hz^2 
-\frac12\int_{Q_T}\hz\nabla\hz\cdot \nabla{\rm div}\vh
\\
& \quad -\int_{Q_T}\Big(\psi\chi\tg(z_t)
+[\![\Delta,M_\psi]\!]z -\psi\tfw(z)\Big)
\big[\vh\cdot\nabla\hz + \frac{\hz}2 ({\rm div}\vh-\rho)\big]
\\
& \quad -\bsC\int_{Q_T}\Big(\phi\chi\tg(z_t)
+[\![\Delta,M_\phi]\!]z -\phi\tfw(z)\Big)\cz
\\
& \quad +\,(BT)_{\Sigma} - {C}_{0,T}^w\,,
\end{split}
\end{equation}
where
\begin{equation} \label{bt-0}
\begin{split}
(BT)_{\Sigma} := &\int_{\Sigma}(\vh\cdot \nabla\hz)\,
\frac{\partial \hz}{\partial\nu}
+ \frac12 \int_{\Sigma}\hz \frac{\partial \hz}{\partial\nu} ({\rm div}\vh-\rho)
\\ & \qquad 
+ \frac12 \int_{\Sigma} \big(\hz_t^2-|\nabla\hz|^2-\lambda\hz^2\big)
(\vh\cdot \nu)
+ \bsC \int_{\Sigma} \cz \frac{\partial \cz}{\partial\nu}\,,
\end{split}
\end{equation}
and
\begin{equation}\label{constant-0-wave}
C_{0,T}^w := \Big[\int_{\Omega}\hz_t \,\vh\cdot \nabla\hz \Big]_0^T 
+ \frac12\,\Big[\int_{\Omega}\hz\hz_t \,({\rm div}\vh-\rho)\Big]_0^T
+ \bsC\Big[\int_{\Omega}\cz_t\cz\Big]_0^T\,.
\end{equation}

\end{lemma}

\subsubsection{Plate component}
Let us now turn to the plate equation.
The abstract equation satisfied by the difference of two 
evolution trajectories reads as follows:
\begin{equation} \label{e:plate-difference}
v_{tt}+A_2 v +\tb(v_t) +\beta\kappa N_0^*Az_t=\tfp(v)\,,
\end{equation}
where $\tb$ and $\tfp$ are defined in \eqref{d:btilde} and 
\eqref{d:f2tilde}, respectively, and 
\[
N_0: H^{s}(\Gam_0) \to H^{s+3/2}(\Om)
\]
is the extension operator 
\[
N_0 \,g = w \quad \Longleftrightarrow\quad
\begin{cases}
(\lam I - \Delta) w = 0 &  \text{in $\Om $}
\\[1pt]
\Dn {w} = 0 & \text{on $\Gam_1$} 
\\[1pt]
\Dn{w} = g &  \text{on $\Gam_0$}
\end{cases}\,.
\]

Temporarily assuming the solution is strong, take the  product 
(in $L_2(0,T;\Gamma_0)$) of equation \eqref{e:plate-difference} with $v_t$, thus obtaining 
the following assertion.

\begin{lemma}[Plate Basic identity]
Suppose Assumptions~\ref{a:sources}, \ref{a:dampings-0} hold.
Let $(z,z_t,v,v_t)$ be the difference of strong  evolution trajectories 
$(h,h_t,u,u_t)$ and $(\zeta,\zeta_t,w,w_t)$. 
Then, for any $T>0$ on has 
\begin{multline} \label{i:plate-basic}
\int_0^T \big(|A_2^{1/2}v|^2+ |v_t|^2\big) = -\int_0^T (\tb(v_t),v)
+ \int_0^T (\tfp(v),v) + 2 \int_0^T |v_t|^2
\\
 \quad +\beta\kappa \int_0^T (N_0^*Az,v_t)
-\beta\kappa\big[(N_0^*A z,v)\big]_0^T-\big[(v_t,v)\big]_0^T\,.
\end{multline}

\end{lemma}
Combining the identities \eqref{i:plate-basic} and 
\eqref{i:wave-basic} establishes the following identity 
for the coupled system.
\begin{lemma}[Basic identity for the composite system]
\label{l:full-basic-identity}
Suppose  the Assumptions \ref{a:geometry}, \ref{a:sources}, \ref{a:dampings-0} hold.
Then, for any positive constant $\bsC$ one has
\begin{equation}\label{e:coupled-basic}
\begin{split} 
\int_{Q_T} (J_{\vh}-\rho I)& \nabla \hz \cdot \nabla \hz 
+ 
\frac{\rho}2\,\int_0^T E_{\hz}(t)+ \bsC\int_0^T E_{\cz}(t)+\int_0^T E_v(t) =
\\
& = \lambda \rho \int_{Q_T}\hz^2  + 2\int_{\Sigma_0} v_t^2
+ 2\,\bsC \int_{Q_T}\cz_t^2 
-\frac12\int_{Q_T}\hz\nabla\hz\cdot \nabla{\rm div}\vh
\\
& \quad -\int_{Q_T}\Big(\psi\chi\tg(z_t)
+[\![\Delta,M_\psi]\!]z -\psi\tfw(z)\Big)
\big[\vh\cdot\nabla\hz + \frac{\hz}2 ({\rm div}\vh-\rho)\big]
\\
& \quad -\bsC\int_{Q_T}\Big(\phi\chi\tg(z_t)
+[\![\Delta,M_\phi]\!]z -\phi\tfw(z)\Big)\cz
\\
& \quad -\int_{\Sigma_0} \tb(v_t)v +\int_{\Sigma_0}\tfp(v)v 
+\beta\kappa \int_{\Sigma_0}zv_t
\\
& \quad +\,(BT)_{\Sigma} - C_{0,T}^w - C_{0,T}^p\,,
\end{split}
\end{equation}
where the boundary terms $(BT)_{\Sigma}$ and the constant $C_{0,T}^w$
are defined in \eqref{bt-0} and \eqref{constant-0-wave} respectively, 
while the constant $C^p_{0,T}$ is defined by 
\begin{equation}
\label{constant-0-plate}
C_{0,T}^p:=-\beta\kappa\Big[(N_0^*Az,v)\Big]_0^T-\Big[(v_t,v)\Big]_0^T\,.
\end{equation}

\end{lemma}

The identity \eqref{e:coupled-basic} is the first step in the proof of existence 
of a global attractor for the dynamical system $(\sY,S(t))$ generated by the PDE
problem \eqref{e:pde-system-0}.
In the following section we show---by careful estimates of all the terms which 
occur in its right hand side---that the above formula eventually yields 
the sought-after pointwise (in time) estimate on the quadratic energy $E_{z,v}$
of the coupled PDE system satisfied by the difference of two trajectories.


\section{Existence of a global attractor}
\label{s:existence-attractor}
This section addresses  existence of a global attractor
for the dynamics generated by the evolutionary problem \eqref{e:pde-system-0}, 
and thus culminates with the proof of Theorem~\ref{thm:attractor-existence}.
Among the key properties satisfied by $(\sY,S(t))$ which will enable us 
to establish the existence of a global attractor, the most challenging 
one is {\em asymptotic smoothness}, whose meaning is recorded in  Definition~\ref{d:asymptotic}.
In turn, this property will eventually follow by showing that the 
compactness criterion recorded in Theorem~\ref{t:khan} can be applied. 

\subsection{Energy inequalities, asymptotic compactness}
Starting from the energy identity in \eqref{e:coupled-basic}, we establish 
a preliminary estimate of the integral over $(0,T)$ of the quadratic energy 
of the system (satisfied by the difference of two trajectories).

\begin{proposition}[Intermediate inequalities: the integral of the quadratic energy]
\label{p:intermediate-estimate}~\break Suppose that the Assumptions~\ref{a:geometry}, \ref{a:sources}, \ref{a:dampings-0} hold.
Let $(h,h_t,u,u_t)$ and $(\zeta,\zeta_t,w,w_t)$ be strong evolution trajectories from  distinct initial data $Y_1=(h_0,h_1,u_0,u_1)$, $Y_2=(\zeta_0,\zeta_1,w_0,w_1)$; 
set $z=h-\zeta$, $v=u-w$.
The following statements hold.
\begin{enumerate}[(i)]
\item
For any $T>0$,  any $\epsilon>0$ there exists positive $T$-independent constants  $C_1,C_2,C_3, C$, and constant
$C_{T,\epsilon}$ dependent on $T$ such that
\begin{equation} \label{e:intermediate-ineq-1}
\begin{split}
& (C_\rho-\epsilon) \int_0^T E_{z,v}(t)\le 
-\int_{Q_T}\psi\chi \tg(z_t)
\big[\vh\cdot\nabla\hz+\frac{\hz}2 ({\rm div}\vh-\rho)\big]
\\
& \quad - C_1 \int_{Q_T} \phi\chi\tg(z_t)\cz 
+ C_2\int_{Q_T} \cz_t^2 
+ \int_{Q_T} \psi\tfw(z) \big(\vh\cdot \nabla\hz\big)
\\
& \quad + C_3
\int_0^T\big[1+(b(w_t),w_t)_{\Gamma_0}+(b(u_t),u_t)_{\Gamma_0}\big]\, 
\|v\|_{2-4\delta,\Gamma_0}
\\
& \quad + \int_0^T\sF(\|\Delta w\|^2,\|\Delta u\|^2)
\,\|\Delta v\|\,\|v\| + \eps \int_{\Sigma_0}v_t^2
\\
& \quad 
+ 
C_{T,\epsilon}\|z\|_{C([0,T],L_2(\Omega))}^2
+ C\big[E(0)+E(T)\big]\,,
\end{split}
\end{equation}
where $\sF(\cdot,\cdot)$ is bounded for bounded values of its arguments.
\item As a consequence of the estimate \eqref{e:intermediate-ineq-1} (after some relabeling of constants),
for any $T$ and any $\epsilon>0$  
\begin{equation} \label{e:intermediate-ineq-2}
\begin{split}
& (C_\rho-\epsilon)\int_0^T E_{z,v}(t) \le  -\int_{Q_T}\psi\chi \tg(z_t)
\big[\vh\cdot\nabla\hz+\frac{\hz}2 ({\rm div}\vh-\rho)\big]
\\
& \quad - C_1 \int_{Q_T} \phi\chi\tg(z_t)\cz 
+ C_2\int_{Q_T} \cz_t^2 
+ \int_{Q_T} \psi\tfw(z) \big(\vh\cdot \nabla\hz\big)
\\
& \quad
+ 
C_{T,\epsilon}\|z\|_{C([0,T],L_2(\Omega))}^2
+ C_{\epsilon,T,\sB}\|A_2^{1/2-\delta}v\|_{C([0,T],L_2(\Gamma_0))}^2
\\[1mm]
& \quad + \eps\int_{\Sigma_0}v_t^2 + C\big[E(0)+E(T)\big]\,.
\end{split}
\end{equation}
\end{enumerate}
\end{proposition}

\begin{proof}
To establish the inequality \eqref{e:intermediate-ineq-1} (and, next, 
\eqref{e:intermediate-ineq-2}), we proceed to estimate all the terms 
on the right-hand side (RHS) of \eqref{e:coupled-basic}. 
In doing so, we will exploit the analysis carried out in \cite{chu-las-tou}
as well as the study performed in \cite{toundykov-na07} for the wave 
equation alone.
Of all the needed calculations only a few are given explicitly; the reader is 
referred to \cite{chu-las-tou}, \cite{toundykov-na07}, or \cite{bu-chu-las},
whenever possible.

1. Since by \eqref{e:positive-hessian} $J_{\vh}-\rho I$ is strictly positive-definite, while
$\intT (E_{\hz}(t) +  E_{\cz}(t)) dt$ is equivalent to the integral $\int_0^T E_z(t)$,
the total quadratic energy of the system $E_{z,v}(t)$ 
readily satisfies 
\begin{equation*}
C_\rho \int_0^T E_{z,v}(t) \le \textrm{LHS of \eqref{e:coupled-basic},}
\end{equation*}
for some positive constant $C_\rho$ (the acronym LHS denotes the `left-hand side').

\medskip

2. Let us turn to the terms in the RHS of \eqref{e:coupled-basic}.
We begin by recalling the following estimates, already utilized 
in \cite{bu-chu-las}: 
\begin{equation} \label{e:estimate-1}
\Big|\int_{\Sigma_0} zv_t\Big| 
= \Big|\int_0^T (N_0^*A z, v_t)_{\Gamma_0}\Big| 
\le \epsilon \int_0^T\|A^{1/2}z\|^2+C_\epsilon\int_0^T|v_t|^2\,,
\end{equation}
\begin{subequations}\label{e:estimate-2}
\begin{align}
\Big|\int_0^T (\tfp(v),v)_{\Gamma_0}\Big|
&\le \int_0^T\sF(\|A_2^{1/2} w\|^2,\|A_2^{1/2} u\|^2)\,\|A_2^{1/2}v\|\,\|v\|
\label{e:estimate-2a}\\
& \underbrace{\le}_{\textrm{when $Y_1, Y_2\in \sB$}}  
\epsilon \int_0^T\|A_2^{1/2}v\|^2+C_{\epsilon,\sB}\int_0^T\|v\|^2\,,
\label{e:estimate-2b}
\end{align}
\end{subequations}
\begin{subequations}\label{e:estimate-3}
\begin{align}
\Big|\int_0^T (\tb(v_t),v)_{\Gamma_0}\Big|
& \le  C \int_0^T\big[1+(b(w_t),w_t)_{\Gamma_0}+(b(u_t),u_t)_{\Gamma_0}\big]\, 
\|A_2^{1/2-\delta}v\|_{\Gamma_0}
\label{e:estimate-3a}\\
& \underbrace{\le}_{\textrm{when $Y_1, Y_2\in \sB$}}
(T+C_\sB) \sup_{[0,T]}\|A_2^{1/2-\delta}v\|_{\Gamma_0}\,.
\label{e:estimate-3b}
\end{align}
\end{subequations}
A few comments are in order. 
The inequality \eqref{e:estimate-1} is straightforward.
The estimate \eqref{e:estimate-2a} (where $\sF(\cdot,\cdot)$ 
denotes a real-valued function which is bounded when its arguments are bounded)
can be verified by elementary computations, while \eqref{e:estimate-2b} 
holds as a consequence of the upper bound for the energy of 
solutions starting in a bounded set $\sB$ (see \eqref{e:bounds}), which gives
\begin{equation*}
|(\tfp(v),v)|\le C_{\sB} \|A_2^{1/2}v\| \|v\|\,.
\end{equation*}
The estimate \eqref{e:estimate-3a} was derived in \cite[Lemma~5.3]{bu-chu-las},
by  using the assumption \eqref{h:pdamping} on the damping 
function $b$, and the Sobolev embedding $\sD(A_2^{1/2-\delta})\subset
H^{2-4\delta}(\Gamma_0)\subset C(\overline{\Gamma_0})$, $0<\delta<1/8$ (see \cite[Lemma~5.3]{bu-chu-las} for more details); \eqref{e:estimate-3a} implies \eqref{e:estimate-3b} when $Y_1, Y_2$
belong to a bounded set $\sB$. In fact, from the energy identity \eqref{e:energy-identity}
it follows
\begin{equation*}
\int_0^T\big[(b(u_t),u_t)_{\Gamma_0}+(b(w_t),w_t)_{\Gamma_0}\big]\le 
\sE_u(0)+\sE_w(0)+2c_1\le C_{\sB}\,.
\end{equation*}
The estimates \eqref{e:estimate-2b} and \eqref{e:estimate-3b} allow us 
to obtain \eqref{e:intermediate-ineq-2} from \eqref{e:intermediate-ineq-1}.

\medskip

3. We now show that the spatial traces---collectively included in $(BT)_{\Sigma}$, as defined by \eqref{bt-0}---are ``almost'' \emph{lower order} terms.
By {\em lower order terms} we mean terms which are finite in topologies below the energy level;
to denote them we will generically utilize the acronym ``{\em lot}''.
The ``almost'' qualifier indicates that these terms can be estimated by means of a combination of 
$\epsilon$-times the quadratic energy (this term will be absorbed by the LHS of the inequality), the plate kinetic energy (which may be expressed in terms of the dampings), plus lower order terms; 
see \eqref{e:final-for-bt} below.

To accomplish this goal, we examine either summand occurring in \eqref{bt-0}.
This analysis parallels that performed in \cite{toundykov-na07} for the (uncoupled) wave equation, though in that case the spatial traces either vanished or
produced at most lower order quantities.
In the present case one needs additionally to take into account the coupling
with the plate equation, which is in fact accomplished through boundary traces. 
The explicit computations below are given for the sake of completeness and
the reader's convenience.

To estimate
\begin{equation*}
I_1:= \int_{\Sigma}(\vh\cdot \nabla\hz)\,\frac{\partial \hz}{\partial\nu}
=\int_0^T\int_{\Gamma}(\vh\cdot \nabla\hz)\,\frac{\partial \hz}{\partial\nu}\,,
\end{equation*}
first observe that  by construction the cutoff function $\psi$ 
vanishes on $\Gamma_0$, while the (Neumann) boundary conditions are homogeneous 
on $\Gamma_1$, and we obtain
\begin{equation} \label{e:normal-derivatives}
\frac{\partial \hz}{\partial\nu}=
\frac{\partial \psi}{\partial\nu}z+\psi\frac{\partial z}{\partial\nu}
= \frac{\partial \psi}{\partial\nu}z\,.
\end{equation}
We now make use of the same argument utilized in \cite[\S 7.3]{toundykov-na07}.
Namely, the gradient $\nabla\hz$ is decomposed into its tangential and normal
(to $\Gamma$) components; the symbol $\nabla_{\Gamma}f(p)$ will denote the 
tangential gradient of $f$ in $p$.
Then, 
\begin{equation}\label{e:i-1}
\begin{split}
I_1 & = \int_0^T\int_{\Gamma}\vh\cdot \big[(\nabla\hz\cdot \nu)\nu+
\nabla_{\Gamma}\hz\big]\,\frac{\partial \psi}{\partial\nu} z
\\[1mm]
&= \int_0^T \int_\Gamma 
\Big[(\vh\cdot\nu)\frac{\partial \hz}{\partial\nu}\frac{\partial \psi}{\partial\nu}z
+ \vh\cdot \big(\nabla_{\Gamma}\psi\,z
+\psi\nabla_{\Gamma}z\big)\frac{\partial \psi}{\partial\nu}z\Big]
\\[1mm]
&=\underbrace{\int_0^T \int_\Gamma 
\Big[(\vh\cdot\nu)\Big(\frac{\partial \psi}{\partial\nu}\Big)^2 z^2
+ \Big[(\vh\cdot \nabla_{\Gamma}\psi)\frac{\partial \psi}{\partial\nu}z^2\Big]}_{I_{11}}
+ \underbrace{\frac12\int_0^T\int_\Gamma \psi\frac{\partial\psi}{\partial\nu}
\big(\vh\cdot \nabla_{\Gamma}z^2\big)}_{I_{12}}\,.
\end{split}
\end{equation}
The first summand $I_{11}$ in \eqref{e:i-1} is readily bounded as follows:
\begin{equation} \label{e:i11-estimate}
|I_{11}|\le C \int_0^T\int_{\Gamma}z^2 = C\int_{\Sigma}z^2
\le \epsilon\int_{Q_T}|\nabla z|^2 
+ C_{\epsilon} \underbrace{\int_{Q_T}|z|^2}_{lot(z)}
\end{equation}
(to complete the estimate, we have used a standard interpolation inequality).
To estimate the integral $I_{12}$, we introduce for simplicity of notation  
\begin{equation*}
\mathbf{G}_\psi:= \psi \frac{\partial \psi}{\partial\nu}\vh\,,
\end{equation*}
and compute 
\begin{equation}\label{e:t12}
K_{12} := \int_\Gamma \mathbf{G}_\psi\cdot \nabla_{\Gamma}z^2 
= \int_\Gamma {\rm div}_{\Gamma}\big(\mathbf{G}_\psi z^2\big)
- \int_\Gamma z^2 {\rm div}\mathbf{G}_\psi = - \int_\Gamma z^2 {\rm div}\mathbf{G}_\psi\,,
\end{equation}
where the step from the second to the last equality invokes the fact that  $\mathbf{G}_\psi$ is compactly supported in a (boundary-less) manifold $\Gam$; this result follows as a corollary of, e.g., \cite[Ch.~8, Theorem~6]{spivak}.
The term  $I_{12}$ (which coincides with $\frac12\int_0^T K_{12}$) can be estimated as done for $I_{11}$ in \eqref{e:i11-estimate}, and therefore 
\begin{equation} \label{e:i1-estimate}
|I_1|\le C \int_{\Sigma}|z|^2\le \epsilon\int_{Q_T}|\nabla z|^2 
+ C_{\epsilon} \int_{Q_T}|z|^2\,.
\end{equation}
 
By \eqref{e:normal-derivatives}, for the integral
\begin{equation*}
I_2:= \int_{\Sigma}\hz \frac{\partial \hz}{\partial\nu} {\rm div}\vh
\end{equation*}
we get
\begin{equation} \label{e:i2-estimate}
|I_2|= \Big|\int_0^T\int_{\Gamma_1}
\psi\frac{\partial \psi}{\partial\nu} z^2 {\rm div}\vh\Big|
\le c\int_{\Sigma_1}|z|^2\,.
\end{equation}
The integral
\begin{equation*}
I_3:= -\frac{\rho}2\int_{\Sigma}\hz \frac{\partial \hz}{\partial\nu}
\end{equation*}
can be treated similarly; therefore  
\begin{equation} \label{e:i3-estimate}
|I_3| \le c\int_{\Sigma_1}|z|^2\,.
\end{equation}
By the properties of the cutoff function $\psi$ and of the vector field $\vh$ 
one has immediately  
\begin{equation}\label{e:i4-estimate}
I_4:= \frac12 \int_{\Sigma} \big(\hz_t^2-|\nabla\hz|^2-\lambda\hz^2\big)(\vh\cdot \nu)
= 0\,.
\end{equation}

It remains to estimate the integral
\begin{equation*}
I_5:= \bsC\int_{\Sigma}\cz \frac{\partial \cz}{\partial\nu}\,.
\end{equation*}
Recall that $\cz=\phi z$ and rewrite $I_5$ accordingly;
next, taking into account the properties of the cutoff function $\phi$
and the boundary conditions
\begin{equation*}
\Dn{z}	= \alpha\kappa v_t \quad \text{on}\quad \Gamma_0,\qquad
\Dn{z}	= 0  \quad\text{on} \quad\Gamma_1,
\end{equation*}
we finally obtain
\[
I_5  = \bsC\int_{\Sigma} \phi z \Big( \Dn{z}
+\phi\frac{\partial z}{\partial\nu}\Big)
= \bsC \int_{\Sigma} \phi\frac{\partial \phi}{\partial\nu}z^2
+ \al\kap\,\bsC\int_{\Sigma_0} zv_t\,.
\]
Thus, using once again interpolation arguments, we get
\begin{equation} \label{e:i5-estimate}
|I_5|\le \epsilon \int_0^T \|A^{1/2}z\|_{\Omega}^2 
+ \epsilon\int_0^T \|v_t\|_{\Gamma_0}^2
+ C_{\epsilon,\bsC}\int_{Q_T}z^2\,.
\end{equation}
Combining the five estimates \eqref{e:i1-estimate}--\eqref{e:i5-estimate},
we conclude 
\begin{equation} \label{e:final-for-bt}
\big|(BT)_{\Sigma}\big|\le \epsilon \int_{Q_T}|\nabla z|^2 + \eps\int_{\Sigma_0}v_t^2 + C_{\epsilon,\bsC}\int_{Q_T}z^2\,.
\end{equation}

4. All time-traces (point-wise energy at $t=0,T$) are dominated by 
\[
C[E(0)+E(T)]\,.
\]

5. For the analysis of the remaining terms the reader is referred to 
\cite{toundykov-na07} and \cite{chu-las-tou}.
We point out explicitly that although the integral
\begin{equation*}
\int_{Q_T}[\![\Delta,M_\psi]\!]z\, \vh\cdot\nabla\hz
\end{equation*}
is equivalent to $\int_{Q_T}|\nabla z|^2$ (that is at full energy level), however, due to the fact that the commutator is supported on the set where $\phi=1$ then, following \cite[\S5.1]{chu-las-tou},  this term can be absorbed by $\bsC\int_0^T E_{\cz}(t)$ by selecting 
$\bsC$ sufficiently large (dependent on $\psi$ and $\mathbf{h}$). More specifically, let
\begin{equation*}
\Omega_{\phi} := \phi^{-1}(\{1\})\cap \Omega\,,
\end{equation*}
and recall that by construction $\psi\equiv 1$ on $\Omega\setminus \Omega_{\phi}$;
this implies, by using \eqref{e:commutator},
\begin{equation*}
[\![\Delta,M_\psi]\!]z= \Delta\,\psi z+2\nabla\psi\cdot\nabla z\equiv 0
\qquad \text{on}\quad \Omega\setminus \Omega_{\phi}\,.
\end{equation*}
The above yields
\begin{equation*}
\begin{split}
\int_{Q_T}[\![\Delta,M_\psi]\!]z \vh\cdot\nabla\hz
&=\int_0^T \int_{\Omega_\phi} (\Delta\psi z+2\nabla\psi\nabla z)\vh\cdot \nabla\hz
\\
&\le C\int_0^T \int_{\Omega_\phi} |\nabla z|^2 
= C\int_0^T \int_{\Omega_\phi} |\nabla \cz|^2 \le C\int_0^T E_{\cz}(t)\,,
\end{split}
\end{equation*}
as desired.

6. We finally observe that the integral
\begin{equation*}
\int_{Q_T}\psi \tfw(z)\vh\cdot \nabla \hz\,,
\end{equation*}
is the most challenging, as it is a ``full energy level'' term, since
$\tfw(z)\sim |z|_{1,\Omega}$.
However, it will not be difficult to cope with this issue at the first stage 
of showing asymptotic smoothness of the semi-flow (thus, existence of the global 
attractor), because of the relatively weak requirements of the abstract result recalled
in Theorem~\ref{t:khan}.
This fact will be clarified below. 
\end{proof}


\medskip
As it will play a fundamental role in the subsequent discussion, let us record the energy 
relation pertaining to the differences of strong trajectories, which we denote 
by $E_{z,v}(t)$:
\begin{multline} \label{e:energy-identity-diff}
E_{z,v}(t) + \beta \int_s^t (\chi \tg(z_t),z_t)\,dr
+ \alpha\int_s^t(\tb(v_t),v_t)\,dr =
\\
E_{z,v}(s) + \beta\int_s^t (\tfw(z),z_t) \, dr
+ \alpha \int_s^t (\tfp(v),v_t) \, dr\,.
\end{multline}
Indeed, the above identity yields the exact expression of the integral of 
the dampings in terms of (pointwise) values of the energy and 
integrals of the nonlinear forces.

\begin{proposition}[Pointwise estimate of the total energy] 
\label{p:energy-pointwise-estimate}
Under the  Assumptions of Proposition~\ref{p:intermediate-estimate}, but now with $(h,h_t,u,u_t)$ and $(\zeta,\zeta_t,w,w_t)$ being generalized solutions (corresponding to initial data $Y_1$ and $Y_2$, respectively) originating in  a bounded set $\sB$, for any sufficiently large $T$ and any $\overline{\epsilon}>0$,
there exists a constant $C_{\sB,\overline{\epsilon}}$ such that 
\begin{equation}\label{e:pointwise-estimate}
E_{z,v}(T)\le C\big(\overline{\epsilon}^2 + 
\max\{g(\pm\overline{\epsilon}),-g(-\overline{\epsilon})\}^2\big)
+ \frac{C_{\sB,\overline{\epsilon}}}T + \Psi_{\sB,T}(Y_1,Y_2)\,,
\end{equation}
where  (with $z= h-\zet$, $v= u-w$)
\begin{multline}\label{d:psi} 
\Psi_{\sB,T}(Y_1,Y_2):=C_T\Big( \sup_{t\in [0,T]}\|z(t)\|^2
+ \sup_{t\in [0,T]}\|A_2^{1/2-\delta}v(t)\|^2
+ \sup_{t\in [0,T]}\|A_2^{1/2-\delta}v(t)\|\Big)
\\
+ C_T \int_0^T ds \int_s^T (\tfw(z),z_t)\, d\tau
+ C_T \int_0^T ds\int_s^T (\tfp(v),v_t) \, d\tau
\\
+ C_T \int_{Q_T} \psi^2 \tfw(z) \vh\cdot \nabla z\,.
\end{multline}
\end{proposition}

\begin{proof}
To derive \eqref{e:pointwise-estimate} from \eqref{e:intermediate-ineq-2}, 
we utilize the chief arguments of \cite[Proof of Lemma~4.2]{chu-las-tou}. First temporarily assume that the trajectories are strong. Begin with analysis of the terms in \eqref{e:intermediate-ineq-2} which involve 
the damping functions.

\textbf{1. (Plate damping.)}
By Assumption~\ref{a:dampings-0}, in particular, according to the lower bound
\eqref{h:pdamping}, it follows that given $\overline{\epsilon}>0$, there exists 
a constant $m_{\overline{\epsilon}}$ such that 
\begin{equation}\label{e:from-hpdamping}
s^2\le \frac1{m_{b,\overline{\epsilon}}} s\,b(s)\,,
\qquad |s|\ge \beps\,.
\end{equation}
Since $v=u-w$, if we set
\begin{equation*}
\Sigma_{0,\overline{\epsilon}}:=\big\{(\xi,t)\in \Gamma_0\times [0,T]:
\; |u_t(\xi,t)|+|w_t(\xi,t)|< \overline{\epsilon}\big\}\,,
\end{equation*}
in light of \eqref{e:from-hpdamping}, we obtain 
\begin{equation*}
\int_{\Sigma_0}v_t^2 \le 2\int_{\Sigma_0}(u_t^2+w_t^2) 
\le 4\int_{\Sigma_{0,\overline{\epsilon}}} \overline{\epsilon}^2 
+ 2m_{b,\overline{\epsilon}}^{-1}\int_{\Sigma_0\setminus \Sigma_{0,\overline{\epsilon}}}
\big(u_tb(u_t)+w_tb(w_t)\big)\,,
\end{equation*}
which implies the estimate 
\begin{equation}\label{e:term-v}
\int_{\Sigma_0}v_t^2\le \overline{\epsilon}^2\,C\,T + C_{\sB,\overline{\epsilon}}\,,
\end{equation}
where $C=4|\Gamma_0|$ is crucially independent of $T$ and $\overline{\epsilon}$. The time-independent bound on the damping integral follows from the energy identity \eqref{e:energy-identity} and the global bound on the energy  from the Proposition \eqref{p:energy-bound}.

\medskip

\textbf{2. (Wave damping.)} The estimate of the three terms involving the wave damping, which occur
on the RHS of \eqref{e:intermediate-ineq-2}, was carried out in \cite{chu-las-tou}.
We provide a few hints for the reader's convenience.
First, notice that 
\begin{multline}\label{e:three-terms}
-\int_{Q_T}\psi\chi \tg(z_t)
\big[\vh\cdot\nabla\hz+\frac{\hz}2 ({\rm div}\vh-\rho)\big]
- C_1\int_{Q_T} \phi\chi\tg(z_t)\cz 
+ C_2\int_{Q_T} \cz_t^2 
\\
\quad \le \epsilon \int_{Q_T}|\nabla z|^2 + C_{T,\epsilon} \sup_{t\in [0,T]}\|z(t)\|^2
+ C_{\epsilon}\int_{Q_T}\big[z_t^2 + \chi \tg(z_t)^2\big]\,.
\end{multline}
To estimate the last integral on the RHS of \eqref{e:three-terms},
we recall that $z=h-\zeta$, and make use of the elementary inequality 
\begin{equation} \label{e:elementary}
s^2+g^2(s) \le \big(M_g+m_{g,\overline{\epsilon}}^{-1}\big)s\,g(s)\,,
\qquad |s|\ge \overline{\epsilon}
\end{equation}
which follows from the Assumption~\ref{a:dampings-0} as well.
Thus, introducing the set
\begin{equation*}
Q_{T,\overline{\epsilon}}:=\{(x,t)\in \Omega_{\chi}\times [0,T]:
\; |h_t(x,t)|+|\zeta_t(x,t)|< \overline{\epsilon}\}\,,
\end{equation*}
and the relative splitting of the integral as before, we thereby obtain
\begin{align}
& \int_{Q_T}\chi \big[z_t^2 + \tg(z_t)^2\big]
\le 2\int_{Q_T}\chi \big(h_t^2+\zeta_t^2 +  g(h_t)^2 + g(\zeta_t)^2\big)
\notag\\
& \quad \le 4\int_{Q_{T,\overline{\epsilon}}} 
\big(\overline{\epsilon}^2 + \max\{g(\overline{\epsilon}),-g(-\overline{\epsilon})\}^2\big)
+ C_{\overline{\epsilon}}\int_{Q_T\setminus Q_{T,\overline{\epsilon}}}
\chi\big(h_tg(h_t)+\zeta_tg(\zeta_t)\big)
\label{e:splitting-w}\\[1mm]
& \quad \le C\,T\big(\overline{\epsilon}^2 
+ \max\{g(\overline{\epsilon}),-g(-\overline{\epsilon})\}^2\big)
+ C_{\overline{\epsilon},\sB}\,,
\label{e:after-splitting}
\end{align}
where, once more, $C=4|\Omega_\chi|$ is independent of $T$ and $\overline{\epsilon}$.

\medskip

\textbf{3. (Energy level wave term.)}
We rewrite  
\begin{equation}\label{e:sum-for-challenging}
\int_{Q_T}\psi \tfw(z)\vh\cdot \nabla \hz=
\underbrace{\int_{Q_T}\psi z\tfw(z)\vh\cdot \nabla \psi}_{I_1}+
\underbrace{\int_{Q_T}\psi^2 \tfw(z)\vh\cdot \nabla z}_{I_2}\,,
\end{equation}
and assert that the first summand satisfies 
\begin{equation}\label{e:final-for-challenging}
\big|I_1\big|\le \epsilon\int_0^TE_z(t) \, dt 
+ C_{\epsilon,\sB,T}\|z\|_{C([0,T],L_2(\Omega))}^2\,.
\end{equation}
(The proof of \eqref{e:final-for-challenging} is relegated to the Appendix.)
The above estimate shows that $I_1$ in \eqref{e:sum-for-challenging} is 
an ``almost lower order'' term (and hence, innocuous): namely,
it is dominated by a term which can be moved to the LHS of 
\eqref{e:intermediate-ineq-2} as well, plus a lower order term.

\medskip

\textbf{4.} We now fix $t=T$ in the energy equality \eqref{e:energy-identity-diff}
(pertaining to the differences of trajectories), and integrate both sides of 
the equality between $0$ and $T$, thus obtaining 
\begin{multline} \label{e:after-integration}
TE_{z,v}(T) \le \int_0^T E_{z,v}(s)\,ds 
\\
+ \beta\int_0^T ds \int_s^T (\tfw(z),z_t)\, d\tau
+ \alpha \int_0^T ds\int_s^T (\tfp(v),v_t) \, d\tau\,.
\end{multline}
Applying all the inequalities \eqref{e:term-v}, \eqref{e:three-terms},
\eqref{e:after-splitting}, (the identity \eqref{e:sum-for-challenging}) and 
\eqref{e:final-for-challenging} to estimate the integral of the quadratic 
energy on the right hand side of \eqref{e:after-integration},
and dividing both sides by $T$, we establish \eqref{e:pointwise-estimate}
for strong solutions.
However, since each term of \eqref{e:pointwise-estimate} is continuous with respect to 
the finite energy topology of $\sY$, the estimate is extended to generalized solutions,
which concludes the proof.
\end{proof}

The existence of a global attractor will follow if we apply the Theorem \ref{t:khan}, as stated in the Appendix. The hypothesis of this theorem directly follows for small 
$\bar{\eps}$ and large enough $T$, provided we also show that the sequential limit 
\eqref{e:sequential-limit} does hold with $\Psi$ as in \eqref{d:psi}.

\begin{proposition}[Weak sequential compactness]
\label{p:compact}
Let $Y^n:=\{\zet^n, \zet_t^n, w^n, w^n_t\}$ be a sequence of trajectories originating  in a bounded subset $\sB$ of $\sY$. Then
\[
\liminf_{n\to\infty}\liminf_{m\to\infty}\Psi_{\sB,T}(Y^n, Y^m) = 0
\]
\end{proposition}
\begin{proof}
First, define
\[
\zmn \dfn \zet^m - \zet^n,\qquad \vmn \dfn w^m - w^n
\]

\noindent\textbf{Step 1: The limits of the lower-order norms.} Recall the following
compactness result (for instance, see \cite{simon}): given a tower of Banach spaces 
$X_0 \stackrel{\mathrm{compact}}{\into}
X\into X_1$, sets which are bounded in $L^p(0,T;X_0)\cap W^{1,r}(0,T;X_1)$
are compact in $L^p(0,T;X)$ for $1\leq p <\infty$ if $r=1$, and $p=\infty$ if $r>1$.
Take
\begin{itemize}
	\item  $X_0= H^1(\Om)$,  $X=L^2(\Om)=X_1$,

	\item and then $X_0 = H^2(\Gam_0)$, $X=\sD(A_2^{1/2-\del})$, 
$X_1=L^2(\Gam_0)$, any $0<\del\leq 1/2$,
\end{itemize}
to conclude that (on a subsequence reindexed again by $n$) $\{\zet^n, w^n\}$ converges strongly in 
$L^\infty(0,T;L^2(\Om)\times \sD(A_2^{1/2-\del}))$ to some $\{\zet, w\}$. In addition, the functions are continuous on $[0,T]$ so
\begin{equation}\label{p:compact:lim:C}
\{\zet^n, w^n\} \to \{\zet, w\}\quad\text{strongly in}\quad 
C([0,T],H^{1-\del}(\Om)\times \sD(A_2^{1/2-\del}))\,.
\end{equation}
Henceforth we will not explicitly mention every passage to a subsequence and continue working with indices labeled $m$ and $n$. Consequently, 
\begin{equation}\label{p:compact:liminf1}
\liminf_{m\to\infty}\liminf_{n\to\infty}\big(\sup_{[0,T]} \|\zmn(t)\| + 
\sup_{[0,T]}\|A_2^{1/2-\del} \vmn\| \big) = 0\,.
\end{equation}

\bigskip

\noindent\textbf{Step 2: Convergence of the source terms.} Pick $t\in[0,T]$, then from the Lipschitz property of $f_1$ we have for every $t$
\[
\begin{split}
\intOm  |f_1(\zet^n(t)) - f_1(\zet(t))| 
\leq & C_{\sB}\, \|\zet^n -\zet\|_{C([0,T],L^2(\Om))}\,.
\end{split}
\]
Hence, by \eqref{p:compact:lim:C} 
\begin{equation}\label{p:compact:lim:f1:C}
 f_1(\zet^n)\to f_1(\zet) \quad \text{strongly}\quad C([0,T],L^1(\Om))\,.
\end{equation}
An almost identical estimate carried with the anti-derivative $F_1$ of $f_1$  shows
\begin{equation}\label{p:compact:lim:F1}
F_1(\zet^n)\to F_1(\zet) \quad \text{strongly}\quad C([0,T],L^1(\Om))\,.
\end{equation}
Let $\chi_s$ be the characteristic function of the set $[s,T]\subset [0,T]$.
Since the sequence $\{\chi_s(t) f_1(\zet^n(x,t))\}$ is bounded in 
$L^\infty([0,T]_s\times [0,T],L^2(\Om))$ (and converges a.e. to $\chi_s f(\zet)$ as follows a fortiori from \eqref{p:compact:lim:f1:C}), then 
\begin{equation}\label{p:compact:lim:f1}
   \chi_s f_1(\zet^n) \to \chi_s f_1(\zet)\quad \text{weakly in }\quad L^2([0,T]_s\times [0,T]\times \Om),
\end{equation}
and trivially
\begin{equation}\label{p:compact:lim:zet-t}
\zet^n_t \to \zet \quad\text{weakly in}\quad L^2([0,T]_s\times [0,T]\times \Om)\,.
\end{equation}

\bigskip

Now establish similar convergence results for the plate component. Since
\begin{equation}\label{p:compact:lim:w:C}
w^n \to w \quad\text{strongly in}\quad C([0,T],H^1(\Gam_0))\,,
\end{equation}
according to \eqref{p:compact:lim:C} with $\del=1/4$, then 
$\|\grad w^n\|^2_{L^2(\Gam_0)} \to \|\grad w\|_{L^2(\Gam_0)}$ in $C([0,T])$. Because a priori  $w^n$ converges weakly in $L^2(0,T;\sD(A_2^{1/2}))$, then $\Del w^n$ converges weakly in $L^2(Q_T)$. As above, we may extend the space-time domain along another dimension with the interval $[0,T]$ to accommodate a characteristic function 
$\chi_s(t)$ of the set $[s,T]\subset [0,T]$; the latter set, being bounded, does not affect the $L^2$ convergence on the finite measure space $[0,T]_s\times[0,T]\times \Om$. Obtain:
\begin{equation}\label{p:compact:lim:w}
\begin{split}
(\chi_s)\left(Q - \|\grad w^n\|^2\right) \Del w^n  \to (\chi_s)\left(Q-\|\grad w\|^2\right)\Del w \\
\text{weakly in}\quad L^2([0,T]_s\times [0,T]\times \Om)
\end{split}
\end{equation}
and 
\begin{equation}\label{p:compact:lim:w-t}
w^n_t \to w \quad\text{weakly in}\quad L^2([0,T]_s\times [0,T]\times \Om)\,.
\end{equation}

\medskip

\noindent\textbf{Step 3: The limits of the source terms.} 
\[
\begin{split}
&\intT ds \int_s^T \big( \tfw(\zmn),\zmn_t\big) dt \\
=&\intOm\intT ds \int_s^T  \frac{d}{dt} \left(F_1(\zet^n) +  F_1(\zet^m)\right) dt-\intT ds \int_s^T  \big[(f_1(\zet^n),\zet^m_t) + (f_1(\zet^m),\zet^n_t) \big]dt\\
=& T \intOm \big[F_1(\zet^n(T)) +F_1(\zet^m(T))\big]dx-\intQT \big[F_1(\zet^n) + F_1(\zet^m)\big]\\
&-\intOm \intT ds \int_s^T \big[ \chi_s f_1(\zet^n)\zet^m_t + \chi_s f_1(\zet^m)\zet^n_t \big].
\end{split}
\]
Pass to the limit $n\to \infty$ and then $m\to \infty$ (on appropriate subsequences),  apply the convergence results \eqref{p:compact:lim:F1}, \eqref{p:compact:lim:f1},  
\eqref{p:compact:lim:zet-t} to obtain on the RHS of the last equality the following terms
\[
2T \intOm F_1(\zet(T))dx- 2\intQT F_1(\zet) -2\intOm \intT ds \int_s^T \frac{d}{dt} F_1(\zet)
\]
which readily cancel each other. Whence
\begin{equation}\label{p:compact:liminf2}
\liminf_{m\to\infty}\liminf_{n\to \infty}
\intT ds \int_s^T \big( \tfw(\zmn),\zmn_t\big) = 0\,.
\end{equation}

Next, recall that $\psi$ (when restricted to the boundary)  is supported on the set where $(\vh\cdot \nu)=0$, consequently $\psi^2 (\vh\cdot \nu)\big|_{\Gam}\equiv 0$, and integration by parts yields
\[
\int_{Q_T}\psi^2 \vh \cdot\grad \big(F_1(\zet^n) + F_1(\zet^m)\big)=
\int_{Q_T}\psi^2 (\div \vh) \big(F_1(\zet^n) + F_1(\zet^m)\big)\,;
\]
use the latter identity to derive:
\[
\begin{split}
\int_{Q_T} \psi^2 \tfw(\zmn) \vh\cdot \nabla \zmn=&
-\int_{Q_T}\psi^2  (\div \vh) \big(F_1(\zet^n) + F_1(\zet^m)\big)\\
&-\int_{Q_T} \psi^2 \vh\cdot \big[f_1(\zet^n)  \nabla \zet^m +
f_1(\zet^m) \nabla \zet^n\big]\,.
\end{split}
\]
Since the sequences $\{\grad \zet^n\}$  and $\{f(\zet^n)\}$ are pre-compact in $L^2(Q_T)$, then passing to the limits $n\to \infty$, $m\to \infty$ in the last identity, and subsequent integration by parts, show
\begin{equation}\label{p:compact:liminf3}
\liminf_{m\to\infty}\liminf_{n\to \infty}\int_{Q_T} \psi^2 \tfw(\zmn) \vh\cdot \nabla \zmn = 0\,.
\end{equation}

\medskip

Finally,
\[
\begin{split}
&\int_0^T ds\int_s^T (\tfp(\vmn),\vmn_t)=\\
=& T\intOm \left[ \frac{Q}{2}\left(\|\grad w^n(T)\|^2 + \|\grad w^m(T)\|\right)^2 - \frac{1}{4}
\left(\|\grad w^n(T)\|^4 +\|\grad w^m(T)\|^4\right)\right]\\
&- \intQT \left[\frac{Q}{2}\left(\|\grad w^n\|^2 + \|\grad w^m\|\right)^2 - \frac{1}{4}
\left(\|\grad w^n\|^4 +\|\grad w^m\|^4\right)\right]\\
&-\intQT ds\int_s^T
\left[\left(Q-\|\grad w^n\|^2\right)\Del w^n w_t^m+
\left(Q-\|\grad w^m\|^2\right)\Del w^m w_t^n\right]\,.
\end{split}
\]
The convergence results \eqref{p:compact:lim:w:C}, \eqref{p:compact:lim:w}, and \eqref{p:compact:lim:w-t} show that the RHS of the last expression converges to
\[
\begin{split}
T\intOm \left[ Q\|\grad w(T)\|^2 - \frac{1}{2}\|\grad w(T)\|^4 \right]
&- \intQT \left[Q \|\grad w\|^2  - \frac{1}{4} \|\grad w\|^4\right]\\
&-2\intT ds\int_s^T \left(Q-\|\grad w\|^2\right)\Del w\, w_t = 0
\end{split}
\]
(use integration by parts in space on the last term to exhibit cancelation). Whence
\begin{equation}\label{p:compact:liminf4}
\liminf_{m\to\infty}\liminf_{n\to \infty}\intT ds \int_s^T (\tfp(\vmn), \vmn_t) = 0\,.
\end{equation}
The limits \eqref{p:compact:liminf1}, \eqref{p:compact:liminf2}, \eqref{p:compact:liminf3} and \eqref{p:compact:liminf4} complete the proof of Proposition \ref{p:compact}.
\end{proof}

\smallskip
The results of the Propositions \ref{p:energy-pointwise-estimate} and \ref{p:compact} confirm, via Theorem \ref{t:khan}, that the dynamical system $(\sY,S(t))$ generated by the PDE \eqref{e:pde-system-0} is asymptotically smooth (see the Definition \ref{d:asymptotic} in the Appendix) and  the existence of a global compact attractor $\mathscr{A}$, along with the claimed geometric description, will follow from the abstract results pertaining  to infinite-dimensional dynamical systems which we summarize below.

\subsection{Concluding the proof of Theorem~\ref{thm:attractor-existence} (existence and geometry of the attractor)}
All the assertions of Theorem~\ref{thm:attractor-existence} will follow 
from \cite[Corollary~2.29]{chu-las-memoirs}. 
To prove that the latter result applies, we need to check that $(\sY,S(t))$ 
possesses three chief properties: namely, that (i) it is {\em gradient}, 
(ii) it is {\em asymptotically smooth} and (iii) the set of its equilibria 
is bounded.
Since the asymptotic smoothness property has been established in the previous 
section, it remains to show that (i) and (iii) hold true.

(i) Let us recall that a dynamical system $(\sY,S(t))$ is gradient if it admits
a strict Lyapunov function. 
We will show that in the present case the Lyapunov function's role is played 
by the full energy $\sE(t)$ of the system.

\smallskip
\noindent
1. We first observe that the identity \eqref{e:energy-identity} 
shows that the map 
\begin{equation*}
t\mapsto \sE(\zeta(t),\zeta_t(t),w(t),w_t(t))
\end{equation*}
is non-increasing in $t$ along strong solutions.
By continuity of $\sE(t)$ in the finite energy norm, this property 
is inherited by weak solutions.

\smallskip
\noindent
2. We further need to show that if 
\begin{equation*}
\sE(\zeta(t),\zeta_t(t),w(t),w_t(t)) = \sE(\zeta(0),\zeta_t(0),w(0),w_t(0))
\qquad  \forall t>0\,,
\end{equation*} 
then $(\zeta(t),\zeta_t(t),w(t),w_t(t))$ is a stationary solution 
of system \eqref{e:pde-system-0}.
Notice, preliminarly, that a stationary solution of the coupled system 
\eqref{e:pde-system-0} has the form $(\zeta,0,w,0)$, where $\zeta$ and $w$ 
satisfy, respectively, the {\em decoupled} boundary value problems 
\begin{equation*}
\begin{cases}
-\Delta \zeta = f(\zeta) & \textrm{in} \;\Omega
\\
\frac{\partial\zeta}{\partial\nu}=0 & \textrm{on} \;\Gamma
\end{cases}\,,
\qquad \qquad
\begin{cases} 
\Delta^2w= f_2(w) & \textrm{in} \;\Gamma_0
\\
w=\Delta w=0 & \textrm{on}\; \partial\Gamma_0
\end{cases}\,.
\end{equation*}
Thus, suppose we are given a generalized solution 
$y(t)=(\zeta(t),\zeta_t(t),w(t),w_t(t))$ such that $\sE_y(t)=\sE_y(0)$ 
for all $t>0$, and let $y^n(t)=(\zeta^n(t),\zeta^n_t(t),w^n(t),w^n_t(t))$ 
any sequence of strong solutions convergent to $y$ in $C([0,T],\sY)$.
From the energy identity it follows that both
\begin{equation*}
\int_0^T\!\int_\Omega\chi g(\zeta^n_t)\zeta^n_t\to 0\,, 
\qquad 
\int_0^T\!\int_{\Gamma_0} b(w^n_t)w^n_t\to 0\,,
\end{equation*}
which implies $\zeta^n_t\to 0$ pointwise a.e. in $\Omega_\chi$ 
and $w^n_t\to 0$ pointwise a.e. in $\Gamma_0$.
That the first limit implies $\zeta_t\equiv 0$ for all positive $t$---i.e. $\zeta$ is 
indeed stationary---follows by applying a unique continuation result established
in \cite{las-trig-zhang};
see also  \cite{toundykov-na07} and \cite[Section~2.7]{chu-las-tou} for further 
references.
\\
More easily, since the sequence $w^n_t$ converges to $w_t$ in 
$C([0,T],L_2(\Gamma_0))$, in view of the second limit we obtain that 
$w_t\equiv 0$ in $\Gamma_0$ for all $t$.
Consequently, $y(t)$ is constant with respect to $t$, i.e. has the form 
$(\zeta,0,w,0)$, that is a stationary solution of \eqref{e:pde-system-0}, 
which concludes the proof.

(iii) Let $\sN$ be the set of equilibria of the flow $S(t)$ associated with 
the PDE problem \eqref{e:pde-system-0}, defined in \eqref{e:equilibria}.
We already observed that $\sN$ is the product of the sets of equilibria 
$\sN_1$, $\sN_2$ of either uncoupled equation.
Boundedness of $\sN_1$ has been shown in \cite[Proposition a-10]{chu-las-tou},
as a consequence of the dissipativity condition in Assumption~\ref{a:sources}.
Boundedness of $\sN_2$ follows as well using the structure of the nonlinear 
function $f_2$ which occurs in the plate equation (according to Berger).
The corresponding proof is fairly simple; however, as it was omitted 
in \cite{bu-chu-las}, it is given here for completeness; for a proof 
under more general assumptions on the nonlinearity, see \cite{potomkin}.

Consider the abstract formulation of the stationary problem corresponding to the plate
equation, that is
\begin{equation*}
A_2w= p_0 + (Q-|A_2^{1/4}w|^2)\, A_2^{1/2} w\,.
\end{equation*}
Taking the dot product (in $L_2(\Gamma_0)$) of this equation by $w$, we obtain 
\begin{equation*}
|A_2^{1/2} w|^2 + |A_2^{1/4}w|^4 - (p_0,w)=Q\, |A_2^{1/4}w|^2\,,
\end{equation*}
which is equivalent to
\begin{equation}  \label{e:to-be-discussed}
|A_2^{1/2} w|^2 + \big(|A_2^{1/4}w|^2-Q\big)\,|A_2^{1/4}w|^2 = (p_0,w)\,.
\end{equation}
Since $A_2 =\Delta^2$ with  hinged boundary conditions is a positive operator, there exists 
$\lambda_0>0$ such that $(A_2w,w)\ge \lambda_0|w|^2$ for any $w$.
Thus, if $Q\le 0$ then \eqref{e:to-be-discussed} readily implies 
the inequality
\begin{equation*} 
\lambda_0|w|^2 + \big(\lambda_0^{1/2}|w|^2-Q\big)\,\lambda_0^{1/2}|w|^2 \le |p_0|\,|w|
\qquad \forall w\,,
\end{equation*}
which cannot hold true, unless there exists a constant $C>0$ such that 
$|w|\le C$.
To show that the same is true when $Q>0$, we proceed by contradiction: 
namely, we assume that for any $C>0$ there exists $w$ such that $|w|>C$.
Thus, taking initially $w$ such that $\lambda_0^{1/2}|w|^2-Q>Q$, and 
next strenghtening the lower bound on $|w|$ to
\begin{equation}\label{e:lower-bound}
|w|> C:=\max\left \{\frac{\sqrt{2Q}}{\lambda_0^{1/4}},\frac{|p_0|}{Q \lambda_0^{1/2}}\right\}\,,
\end{equation}
we see that \eqref{e:to-be-discussed} combined with \eqref{e:lower-bound} imply
\begin{equation*}
0<\lambda_0|w|^2 + |w|\big(Q\lambda_0^{1/2}|w|-|p_0|\big)\le 0\,,
\end{equation*}
and we achieve a contradiction. 
\qed


\section{Finite dimensionality and smoothness of the attractor}
\label{s:finite-dim}
In this section we discuss the issue of fractal dimension of the global attractor.
We aim to show that the criterion recalled as Theorem~\ref{t:ladyzhenskaia-general}
applies, i.e.~the global attractor has finite fractal dimension,
thereby establishing Theorem~\ref{thm:finite-dim}.
As we will see, while the  ``Carleman version'' of the wave fundamental 
identity \eqref{e:fundamental-carleman} has not been used in order to 
establish the existence of the attractor, it is central to the proof of finite dimensionality of the attractor.
Indeed, it brings about the following Lemma, which constitutes a major 
step in the proof of the inequality \eqref{e:involves-seminorms} required by
Theorem~\ref{t:ladyzhenskaia-general}.

\begin{lemma}[Observability/stabilizability-like inequalities] 
\label{l:central-to-finitedim}
Suppose the Assumptions~\ref{a:geometry}, \ref{a:sources} hold, along with,
initially, the weaker Assumption~\ref{a:dampings-0} on the damping 
functions.  Let $Y_1=(h_0,h_1,u_0,u_1)$, $Y_2=(\zeta_0,\zeta_1,w_0,w_1)$, be two smooth initial data and the  $(h,h_t,u,u_t)$ and $(\zeta,\zeta_t,w,w_t)$ be the 
corresponding strong evolution trajectories; set $z=h-\zeta$, $v=u-w$.
(Recall the definitions of the cutoff functions $\psi$ and $\phi$, the vector
field $\vh$, as well as $\hz=\psi z$ and $\cz=\phi z$.)

\textsc{Then} for a sufficiently large $T$ and any positive parameters $\tau$, $\epsilon$, 
$\epsilon_0$, there exist positive constants $C_1$, $C_2$, $C_{1;T}$, $C_{2;T}$, and  
$C_{\tau,T,\epsilon,\epsilon_0}$ such that 
\begin{equation} \label{e:observability-0}
\begin{split}
& C_1\bigg\{\int_{Q_T} e^{\tau\Phi}\big(|\nabla \hz|^2+ \lambda \hz^2 + \hz_t^2\big)
+\int_{Q_T} e^{\tau\Phi}\big(|\nabla \cz|^2 +\lambda \cz^2 + \cz_t^2\big)\bigg\}
\\
& \qquad \le \,\frac{\epsilon_0^2}4\int_0^T\|\nabla z\|^2 
+\epsilon\int_{Q_T} e^{\tau\Phi}\big(|\nabla z|^2+ z_t^2+\tfw(z)^2\big)
\\
& \qquad\qquad -\tau \int_{Q_T}e^{-\tau\Phi}\sM_1^2 
+ \int_{Q_T}\big[\psi \tfw(z)-\chi \tg(z_t)\big]\,\sM_1
\\
& \qquad\qquad+ C_2 \int_{Q_T}e^{\tau\Phi}\chi^2\tg(z_t)^2 
+C_{1;T}\int_0^T\int_{\Omega_\chi}e^{\tau\Phi}z_t^2
+ \eps\int_{\Sigma_0} v_t^2
\\ 
& \qquad\qquad + C_{2;T}e^{-\delta\tau}\big[E_z(0)+E_z(T)\big]
+ C_{\tau,T,\epsilon,\epsilon_0}\|z\|_{C([0,T],L_2(\Omega))}^2\,.
\end{split}
\end{equation}
\end{lemma}
 
The estimate \eqref{e:observability-0} is an analogue of the
inequality (46) in the first part of \cite[Lemma~4.3]{chu-las-tou}
(which pertained to the wave equation alone).
It is established as well by carefully estimating all the terms which 
occurr in the RHS of the wave fundamental identity \eqref{e:fundamental-carleman}.
Here, due to the coupling, the wave (non-homogeneous) boundary traces naturally
bring in the inequality the integral of the plate kinetic energy. 
The proof is omitted; for technical details see \cite[Section~6.6]{chu-las-tou}, 
and step 3. in the proof of Proposition~\ref{p:intermediate-estimate} of the present paper.

\begin{remark}
\begin{rm}
We just notice that in order to obtain the estimate \eqref{e:observability-0}, the constant $c$ (which occurrs in the definition of the function $\Phi(x,t)$) has been 
assumed to satisfy both the inequalities $c< \rho/2$ and \eqref{e:constraint-2}. 
Let us recall from \cite{chu-las-tou} that in particular, in light of the latter 
constraint one has $\Phi(x,0)=\Phi(x,T)<0$, which in turn implies the existence of a constant $\delta>0$ such that 
\begin{equation}
\Phi(x,0) <-\delta\,, \;\Phi(x,T) <-\delta \qquad \forall x\in \Omega\,.
\end{equation}
The above property has been critically used to estimate the integrals (evaluated
at the end points $t=0,T$) \eqref{e:time-traces} in \eqref{e:fundamental-carleman}.
\end{rm}
\end{remark}

\smallskip

The remainder of the proof is split into a sequence of three results:
\begin{enumerate}
	\item  Lemma~\ref{l:final} provides the inequality which leads to the finite-dimensionality estimate and  the regularity of the attractor; however, the result possibly holds only on a restricted time-interval and under an assumption that a certain auxiliary estimate is true.

	\item Proposition~\ref{p:pre-final} verifies the auxiliary inequality used in the hypothesis of Lemma~\ref{l:final}, but up to a perturbation by source terms. 
	
	\item Finally, Proposition~\ref{p:forces} uses the perturbed estimate derived in  Proposition~\ref{p:pre-final} to show that the hypothesis of Lemma \ref{l:final} holds 
and the necessary estimate is satisfied globally in time, independently of the choice of the trajectories through the attractor.
  
\end{enumerate}

\begin{lemma}[Conditional regularity and observability]\label{l:final}
Let $T$ be given by  Lemma \ref{l:central-to-finitedim}. Suppose  for any pair of trajectories $\gam_1\dfn\{h,h_t, u,u_t\}$ and $\gam_2\dfn\{\zet,\zet_t, w,w_t\}$ through the attractor $\mathscr{A}$ one can find a time $T_{\gam_1,\gam_2}$ and a non-negative function $G \in L^1(\mathbb{R})$ such that
\begin{equation}\label{l:final:assum}
E_{z,v}(s+T)\le \sigma E_{z,v}(s)+C_{\mathscr{A},\sigma,T}lot_s^{s+T}(z,v) + 
\int\limits_s^T G(\tau) E_{z,v}(\tau)d\tau\quad 
\forall s< T_{\gam_1,\gam_2}-T\,,
\end{equation}
where 
\begin{itemize}

\item
$\displaystyle lot_a^b(z,v) \dfn \sup_{\theta\in [a,b]}
\Big(\|z(\th)\|^2 + \|A_2^{1/2-\delta}v(\theta)\|^2\Big)$,

\item the constants $0<\sigma < 1$ and $C_{\mathscr{A},T,\sig}>0$ are \textbf{independent} of $\gam_1, \gam_2$ (in the attractor),

\item $\|G\|_{L^1(\bbR)} <\infty$ and the norm can be bounded \textbf{independently} of 
$\gam_1, \gam_2$.
\end{itemize}

\smallskip

\noindent Then
\begin{enumerate}[(a)]
	\item There exists $\sig_1<1$, $\barT\geq T$ and $\barT_{\gam_1,\gam_2}< T_{\gam_1,\gam_2}$ such that the following estimate holds for any pair of trajectories through the attractor with all the coefficients being \textbf{independent of the  trajectories} themselves:
	\begin{equation}\label{e:final}
	E_{z,v}(s+\barT)\le \sigma_1 E_{z,v}(s)
+C_{\mathscr{A},\sigma,\barT}lot_s^{s+\barT}(z,v),
\quad \forall s< \barT_{\gam_1,\gam_2}\,.
	\end{equation}

	\item Every trajectory $\gam=\{\zet,\zet_t,w,w_t\}$ through the attractor is strong. Furthermore there exists a constant $C_{\mathscr{A}}$, dependent  on the diameter of the attractor $\mathscr{A}$ (in the state space $\sY$), but independent of $\gam$ such that 
\[
\|\zet_t\|_{1,\Om} + \|\zet\|_{2,\Om} + \|w_t\|_{2,\Om} + \|w\|_{4,\Om} 
\leq C_{\mathscr{A}},\quad \forall t < \barT_{\gam}
\]
for some $\barT_\gam$. Moreover, if in part (a) $\barT_{\gam_1,\gam_2}  = \infty$ for every $\gam_1, \gam_2$, then $\barT_\gam = \infty$, i.e. the attractor is bounded in the higher-energy space.
\end{enumerate}
\end{lemma}

\begin{proof}
Let $T$ be given by Lemma \ref{l:central-to-finitedim}, and fix any $T_2\in[T,2T]$ (in fact any interval with $T$ being the left end-point would serve the purpose). Repeatedly applying the inequality \eqref{l:final:assum} from the hypothesis get
\[
\begin{split}
& E_{z,v}(s+T_2) \leq \\
\leq & \sigma E_{z,v}(s)+C_{\scrA,\sigma,T_2}lot_s^{s+T_2}(z,v) + \int_s^{T_2} G(\tau) E_{z,v}(\tau)d\tau \\
\leq  &
\sigma^2 E_{z,v}(s-T_2)+(\sigma +1)\left( C_{\scrA,\sigma,T_2}lot_{s-T_2}^{s+T_2}(z,v) + \int_{s-T_2}^{s+T_2} G(\tau) E_{z,v}(\tau)d\tau\right) \leq \cdots\\
\leq & \sig^m E_{z,v}(s-(m-1)T_2)\\
& + \sum_{j=0}^{m-1} \sig^j\!
\left(C_{\scrA,\sigma,2T}lot_{s-(m-1)T_2}^{s+T_2}(z,v) + \int\limits_{s-(m-1)T_2}^{s+T_2} G(\tau) E_{z,v}(\tau)d\tau\right)\,,
\end{split}
\]
where the very last step also uses the fact that the constant $C_{\mathscr{A},\sig,T_2}$ is continuous increasing with respect to $T_2$, hence can be bounded by 
$C_{\mathscr{A},\sig,2T}$. Define $s_0 \dfn s-(m-1) T_2$ and $t = s_0+mT_2= s+T_2$, then
\[
E_{z,v}(t)\leq \sig^m E_{z,v}(s_0) + 2 C_{\scrA,\sig,2\cdot T} lot_{s_0}^{t}(z,v)
+2 \int_{s_0}^t G(\tau) E_{z,v}(\tau)d\tau.
\]
The latter inequality  depends on $T_2$ only via $t$, hence holds for all $s+ T \leq t\leq s+2T$, or, equivalently,
$s_0+mT\leq t \leq s_0+2mT$. Since $lot_{s_0}^{t}(z,v)$ is non-decreasing in $t$, Gronwall's inequality on this time-interval  interval implies
\[
E_{z,v}(t)\leq \big(\sig^m E_{z,v}(s_0) + 2C_{\mathscr{A},\sig,T}lot_{s_0}^t(z,v)\big)
e^{2\|G\|_{L^1(\bbR)}}.
\]
Because $\sig<1$, we can carry out this argument for $m$ large enough so that 
\[
  \sig_1\dfn \sigma^m \cdot e^{2\|G\|_{L^1(\bbR)}} < 1.
\]
For such an $m$ it then follows
\[
E_{z,v}(t) \leq \sig_1 E_{z,v}(s_0) + C_{\mathscr{A},\sig_1,T} lot_{s_0}^t(z,v),\quad \forall t\in [s_0+mT, s_0+2mT]\,.
\]
To obtain part (a) of the Lemma just relabel the constants. For instance, pick $\barT\in [s_0+mT, s_0+2mT]$, denote $\barT_{\gam_1,\gam_2}\dfn T_{\gam_1,\gam_2}-2mT$, and, finally, relabel $s_0$ into $s$.

\medskip

Now pick any trajectory $t\mapsto\{\zet(t),\zet_t(t),w(t),w_t(t)\}$ through $\mathscr{A}$, let $h\in(0,1)$, and introduce the difference quotients
\[
z^h(t) \dfn \zet(t+h)-\zet(t), \qquad v^h(t) \dfn w(t+h)-w(t)\,.
\] 
According to the now-verified part (a) of the Proposition in question we may find $\barT>0$  and some $\barT_{\zet,w}$ (slightly decreased if needed to accommodate a shift by $h<1$ along the trajectory), and $\sig_1<1$ such that 
\[
E_{z^h,v^h}(s+\barT)  \le \sigma_1 E_{z^h,v^h}(s)  +C_{\mathscr{A},\sigma_1,T}\, 
lot_s^{s+T}(z^h,v^h)\qquad \forall s<\barT_{\zet,w}\,.
\]
Divide now each side of the equation by $h^2$ and introduce
\[
 Y^h(t) \dfn \frac{1}{h}\big\{z^h(t), z^h_t(t), v^h(t), v^h_t(t)\big\}\,,
\]
with the respective energy denoted $E^h(t)$ for a shorthand.  
Observe that
\[
\frac{1}{h^2}lot_a^b(z^h, v^h) 
= \sup_{\th\in[a,b]}\left(
\left\|\frac{\zet(\th+h)-\zet(\th)}{h}\right\|^2 
+ \left\|\frac{w(\th+h)-w(\th)}{h}\right\|^2 \right)\,.
\]
Since both $\zet$ and $w$ are in $C^1([0,T];L^2)$ (respectively over $\Om$ and $\Gam_0$), these difference quotients can be bounded via $\|\zet_t\|^2$ and $\|w_t\|^2$ which, in turn, are uniformly globally bounded by some $C_{\mathscr{A}}$. Consequently, we may without loss of generality state
\[
E^h(s+\barT) \le \sigma_1 E^h(s) +C_{\mathscr{A},\sigma_1,\barT} 
\quad \forall s \leq \barT_{\zet,w},\quad h\in (0,1)\,.
\]	
Both sides of the last inequality are continuous with respect to $s$; take $\sup_{s\leq \barT_{\zet,w}}$ on the RHS, and $\sup_{s\leq  \barT_{\zet,w}-\barT}$ on the LHS, to conclude that 
\[
(1-\sig_1) \sup_{t\leq \barT_{\zet,w}} E^h(t) \leq C_{\mathscr{A},\sig_1, \barT}\,;
\]
this last estimate is independent of $h$, hence taking $h\searrow 0$ yields 
\[
 \|\zet_{tt}(t)\| + \|\grad \zet_t(t)\| + \|w_{tt}(t)\|  + \|\Del w_t(t)\|
 < C_{1,\mathscr{A}}\quad \forall t\leq \barT_{\zet,w}\,.
\]
From the system \eqref{e:pde-system-0} it then follows that the $H^2$ and $H^4$ norms of 
$\zet$ and $w$, respectively, are also bounded by some constant $C_{\mathscr{A}}$, at least for  $t\leq \barT_{\zet,w}$. Forward propagation of regularity implies that the trajectory is strong, in particular that 
\[
\mathscr{A} \subset H^2(\Om)\times H^1(\Om) \times H^4(\Om)\times H^2(\Om)\,;
\]
however, we cannot yet claim that the regularity is uniform since for each trajectory, the bound $C_{1,\mathscr{A}}$ in the higher topology, albeit not directly dependent on $t$, came from the analysis carried out only for until a certain time $t$. 
However, if the original $\barT_{\zet,w}$ provided by part (a) of the Lemma is infinite then taking $\barT_{\zet,w}\to \infty$ implies the said bound for all $t\in \bbR$.
This completes the proof of Lemma \ref{l:final}.
\end{proof}

\bigskip

Now to complete the proof of finite-dimensionality of the attractor and its regularity it remains to verify that the hypothesis of  Lemma \ref{l:final} holds, with $T_{\gam_1,\gam_2}$ of \eqref{l:final:assum} being $\infty$. As a first step we verify the desired  estimate, but perturbed by the source terms.

\begin{proposition}[Perturbed estimate]\label{p:pre-final}
Suppose the Assumptions~\ref{a:geometry}, \ref{a:sources} hold. Furthermore, let the damping functions $g$ and $b$ further satisfy Assumption \ref{a:dampings-1}, then there exist constants  $0<\sigma<1$, $c_E>0$, $\sC_{\sB,\sigma}$, and $T_0 >0$, such that for any  $T\geq T_0$ and  any two trajectories $Y_1=\{h,h_t, u,u_t\}$, $Y_2=\{\zet,\zet_t, w,w_t\}$  (with $z:= h-\zet$, $v:=u-w$)
\begin{equation}\label{e:pre-final}
E_{z,v}(s+T) +  c_E\intsT E(t) dt\le \sigma E_{z,v}(s) 
+\sC_{\sB,\sigma}\, lot(z,v) + \{\textsc{Forces}\}\quad \forall s\in\bbR\,,
\end{equation} 
where, 
\begin{subequations}\label{d:forces}
\begin{eqnarray}
\textsc{Forces}&\dfn&   \al \{\sI\}(\tfw(z),z_t) + \bet\{\sI\} (\tfp(v),v_t)\\
\{\sI\} &\dfn&\Big\{ C_{1,T}\intsT dt \int_s^t + C_{2,T}\int_{t_0}^{t_1}dt \int_t^{s+T} + C_{3,T}\intsT\Big\}
\end{eqnarray}
\end{subequations}
and $(t_0, t_1)\subset [s,s+T]$.
\end{proposition} 

\begin{remark}
\begin{rm}
The constant $c_E$ in \eqref{e:pre-final} was granted a special name merely to make it easier to keep track of this constant. This parameter will determine how constants are chosen later on in Step 1 of the proof of Proposition \ref{p:forces}.
\end{rm}
\end{remark}

\begin{proof}
1. First work with strong trajectories. For any $\eps>0$ we have
\begin{multline}
-\tau \int_{Q_T}e^{-\tau\Phi}\sM_1^2 
+ \int_{Q_T}\big[\psi \tfw(z)-\chi \tg(z_t)\big]\,\sM_1
\\
\le -\Big(\tau-\frac1{2\epsilon}\Big)\int_{Q_T}e^{-\tau\Phi}\sM_1^2 
+ \epsilon\int_{Q_T}e^{\tau\Phi}\big[\tfw(z)^2+\chi^2\tg(z_t)^2\big]\,.
\end{multline}
Applying this estimate to \eqref{e:observability-0} gives
\begin{equation} \label{e:observability-1}
\begin{split}
& C_1\left\{\int_{Q_T} e^{\tau\Phi}\big(|\nabla \hz|^2+ \hz_t^2\big)
+\int_{Q_T} e^{\tau\Phi}\big(|\nabla \cz|^2+ \cz_t^2\big) \right\}
\\
& \qquad \le \,\frac{\epsilon_0^2}4\int_0^T\|\nabla z\|^2 
+2\epsilon\int_{Q_T} e^{\tau\Phi}\big(|\nabla z|^2+ z_t^2+\tfw(z)^2\big)
-\Big(\tau-\frac1{2\epsilon}\Big)\int_{Q_T}e^{-\tau\Phi}\sM_1^2 
\\
& \qquad\qquad+ C_2 \int_{Q_T}e^{\tau\Phi}\chi^2\tg(z_t)^2 
+C_{1;T}\int_0^T\int_{\Omega_\chi}e^{\tau\Phi}z_t^2
+ \eps\int_{\Sigma_0} v_t^2
\\ 
& \qquad\qquad + C_{2;T}e^{-\delta\tau}\big[E_z(0)+E_z(T)\big]
+ C_{\tau,T,\epsilon,\epsilon_0}\|z\|_{C([0,T],L_2(\Omega))}^2\,.
\end{split}
\end{equation}
Next choose 
$\epsilon=\tau^{-1}$, which ensures
\begin{equation*}
\Big(\tau-\frac1{2\epsilon}\Big)\int_{Q_T}e^{-\tau\Phi}\sM_1^2  > 0\,.
\end{equation*}

\smallskip
\noindent
2.  By construction of $\Phi$  (see Section \ref{s:vec-field}) we have
\[
\Phi(x,0) = \Phi(x,T) < 0 \quad \text{and}\quad 
\Phi(x,T/2) = d(x) \ge \min_{x\in \bar{\Om}}d(x) > 0\,.
\]
So there exists 
an interval $[t_0,t_1]\subset (0,T)$ such that 
\begin{equation}\label{e:positive-weight}
\Phi(x,t) \ge 0 \qquad (x,t)\in \Omega\times [t_0,t_1]\,. 
\end{equation}
Let us choose, specifically,
\begin{equation*}
t_0= \frac{T}2-\sqrt{\frac{\inf_{x\in \overline{\Omega}} d(x)}c}\,, \quad
t_1= \frac{T}2+\sqrt{\frac{\inf_{x\in \overline{\Omega}} d(x)}c}\,. 
\end{equation*}
The above observation enables us to separate the terms $z_t^2$ and 
$|\nabla z|^2$ from the Carleman weight $e^{\tau\Phi}$, thus 
recovering the (wave) energy integrals.
In fact, from \eqref{e:positive-weight} it follows 
$e^{\tau \Phi} \ge 1$ in $\Omega\times [t_0,t_1]$,
and therefore 
\begin{equation*}
C_1 \int_{t_0}^{t_1}\!\!\int_{\Omega}\big(|\nabla z|^2+ z_t^2\big)
\le C_1 \int_{Q_T}e^{\tau\Phi}\big(|\nabla z|^2+ z_t^2\big)\,,
\end{equation*}
which in turn yields 
\begin{equation}\label{e:observ:energy-left}
C_{1,\Omega} \int_{t_0}^{t_1}E_z(t) \, dt
\le C_1 \int_{Q_T}e^{\tau\Phi}\big(|\nabla z|^2+ \lambda z^2+ z_t^2\big)\,.
\end{equation}

\smallskip
\noindent
3. We now manage the integrals involving the dampings.
The stronger Assumption~\ref{a:dampings-1} on the wave damping implies
\begin{equation*} 
s^2+\tg^2(s) \le \big(M'_g+{m'_g}^{-1}\big)s\,\tg(s) \qquad \forall s\,,
\end{equation*}
which enables us to obtain 
\begin{multline} \label{e:weight-free}
C_2 \int_{Q_T}e^{\tau\Phi}\chi^2\tg(z_t)^2 
+C_{2;T}\intT \int_{\Omega_\chi}e^{\tau\Phi}z_t^2
\\
\le C_T \int_{Q_T}e^{\tau\Phi}\chi\tg(z_t)\,z_t
\le C_{\tau,T}\int_0^T (\chi\tg(z_t),z_t)_\Omega\,.
\end{multline}
As for the plate damping, by the lower bound in \eqref{e:plate-hypo-stronger}
we have that  
\begin{equation*} 
s^2 \le {m'_b}^{-1}s\tb(s)\,, \qquad \forall s\,,
\end{equation*}
which immediately gives
\begin{equation}\label{e:kinetic-to-damping}
\int_{\Sigma_0} v_t^2 \le \frac1{m'_b}\int_{\Sigma_0} v_t\,\tb(v_t )\,.
\end{equation}
Use \eqref{e:observ:energy-left}, \eqref{e:weight-free}, \eqref{e:kinetic-to-damping} and choose $\eps$ sufficiently small in \eqref{e:observability-1}, to finally obtain the following 
inequality for the integral (over $(t_0,t_1)$) of the wave quadratic energy:
\begin{multline} \label{e:wave-est}
\int_{t_0}^{t_1}\beta E_z(t)\le \beta\frac{\epsilon_0^2}8\int_0^T\|\nabla z\|^2\,dt 
+ C_{\tau,T,\epsilon_0}\Big\{\underbrace{\int_0^T \big[\beta(\chi\tg(z_t),z_t)_\Omega
+\alpha(\tb(v_t),v_t)_{\Gamma_0}\big]}_{\textsc{Dampings}}\Big\}
\\
+ C_T\beta e^{-\delta\tau}\big[E_z(0)+E_z(T)\big]
+ C_{\tau,T,\epsilon,\epsilon_0}\|z\|_{C([0,T],L_2(\Omega))}^2\,,
\end{multline}
where $C_T$ is independent of $\tau$. 

\smallskip
\noindent
4. On the other hand, using similar arguments as in \cite[Lemma~4.2]{bu-chu-las}, 
we get the following estimate of the integral (over $(0,T)$) of the 
plate quadratic energy: 
\begin{multline} \label{e:plate-est}
(1-\sigma_2)\int_0^T \alpha E_v(t)\le 
\sigma_1\frac{\beta}2\int_0^T\|\nabla z\|^2 + C_{\sigma_1,\eps}\int_0^T\|v_t\|^2
+ \eps\big[E(0)+E(T)\big]
\\
+C_{\sB,T,\sigma_2,\eps}\|A_2^{1/2-\delta}v\|_{C([0,T],L_2(\Gamma_0))}^2
+C_{\eps,\sB,T} \|z\|^2_{C([0,T],L^2(\Om))}\,.
\end{multline}

\smallskip
\noindent
5. Rewrite the energy relation \eqref{e:energy-identity-diff}
with $s=0$, integrate on $(0,T)$ and multiply both sides by $\epsilon_0$,
thus obtaining ($E(t)$ here stands for  $E_{z,v}(t)$)
\begin{multline*} 
\epsilon_0\int_0^T E(t)\, dt 
+ \epsilon_0\int_0^T\beta \int_0^t (\chi \tg(z_t),z_t)_{\Omega}\,dr\, dt
+ \epsilon_0\int_0^T\alpha\int_0^t(\tb(v_t),v_t)_{\Gamma_0}\,dr\, dt =
\\
= \epsilon_0TE(0) + \epsilon_0\Big\{\underbrace{\int_0^T\beta\int_0^t (\tfw(z),z_t)\,dr\, dt
+ \int_0^T\alpha \int_0^t (\tfp(v),v_t)\,dr\,dt}_{\textsc{Forces I}}\Big\}\,.
\end{multline*}
which readily implies 
\begin{equation} \label{e:sources1}
\epsilon_0\int_0^T E(t)\,dt\le \epsilon_0TE(0)+\epsilon_0 \big\{\textsc{Forces I}\big\}\,.
\end{equation}
Summing \eqref{e:wave-est}, \eqref{e:plate-est} (specifically, 
with $\sigma_1=\epsilon/8$, $\sigma_2=1/2$) and \eqref{e:sources1}, get
\begin{align}\label{e:coupled-est-1} 
& \int_{t_0}^{t_1}\beta E_z(t)+\int_0^T \alpha E_v(t)+\epsilon_0\int_0^T E(t)
\notag\\
& \qquad \le C_{\tau,T,\epsilon,\epsilon_0}\big\{\textsc{Dampings}\big\}
+ C_T\beta e^{-\delta\tau}\big[E_z(0)+E_z(T)\big]
\notag\\
& \qquad\qquad + \underbrace{\big(\frac{\epsilon_0}4+\frac{\epsilon_0}4\big)\int_0^T E(t)}_{\textrm{it can be absorbed by the LHS}} 
+ C_{\epsilon_0}\int_0^T\|v_t\|^2
\notag\\
& \qquad\qquad + 2\eps\big[E(0)+E(T)\big]+ \epsilon_0TE(0) 
+ \epsilon_0 \big\{\textsc{Forces I}\big\}
\notag\\
& \qquad\qquad + C_{\tau,T,\epsilon_0,\eps}\|z\|_{C([0,T],L_2(\Omega))}^2
+C_{\sB,T,\eps}\|A_2^{1/2-\delta}v\|_{C(0,T],L_2(\Gamma_0))}^2\,.
\end{align}
Next, move the $\epsilon_0/2$ integral of the total energy of the system 
to the LHS of \eqref{e:coupled-est-1}, multiply the obtained inequality
by $2$, and add 
\begin{equation*}
2\int_{t_0}^{t_1}\alpha E_v(t)
\end{equation*}
to both sides, thereby obtaining 
\begin{align}\label{e:coupled-est-2} 
& \int_{t_0}^{t_1}E(t)+\epsilon_0\int_0^T E(t)
\le C_{\tau,T,\epsilon,\epsilon_0}\big\{\textsc{Dampings}\big\}
+ C_T\beta e^{-\delta\tau}\big[E_z(0)+E_z(T)\big]
\notag\\
& \qquad + C_{\epsilon_0}\int_0^T\|v_t\|^2
+ 2\eps\big[E(0)+E(T)\big]+ \epsilon_0TE(0) + \epsilon_0 \big\{\textsc{Forces I}\big\}
\notag\\
& \qquad + \underbrace{C_{\tau,T,\eps,\epsilon_0}\|z\|_{C([0,T],L_2(\Omega))}^2
+C_{\sB,T,\eps}\|A_2^{1/2-\delta}v\|_{C(0,T],L_2(\Gamma_0))}^2}_{lot(z,v)}\,.
\end{align}

On the other hand, setting now $t=T$ in the energy indentity 
\eqref{e:energy-identity-diff}
and integrating both sides in $s\in (t_0,t_1)$, we have also
\begin{equation} \label{e:sources2}
(t_1-t_0) E(T) \le \int_{t_0}^{t_1} E(s)\,ds+
\Big\{\underbrace{\int_{t_0}^{t_1}\Big[\beta\int_s^T (\tfw(z),z_t)\,dr
+ \alpha \int_s^T (\tfp(v),v_t)\,dr\Big]\,ds}_{\textsc{Forces II}}\Big\}\,.
\end{equation}
Adding together \eqref{e:coupled-est-2} with \eqref{e:sources2} yields
\begin{multline} \label{e:coupled-est-3} 
(t_1-t_0) E(T) +  \epsilon_0\int_0^T E(t)\, dt \le 
C_{\tau,T,\eps,\epsilon_0}\big\{\textsc{Dampings}\big\}
+ C_T\beta e^{-\delta\tau}\big[E_z(0)+E_z(T)\big]
\\
+ 2 \epsilon_0 TE(0) + 2 \epsilon_0 \big\{\textsc{Forces I}\big\}+
\big\{\textsc{Forces II}\big\}
+ 2\eps \big[E(0)+E(T)\big] + lot(z,v)\,.
\end{multline}

\smallskip
\noindent
6. Next, use a by now standard argument.
Rewrite once again the identity \eqref{e:energy-identity-diff}, this time with
$s=0$ and $t=T$, resulting in an exact expression of the integrals involving 
the dampings:
\begin{multline*}
\beta \int_0^T (\chi \tg(z_t),z_t)\,dt
+ \alpha\int_0^T(\tb(v_t),v_t)\,dt
\\
=E(0)- E(T) + \underbrace{\beta\int_0^T (\tfw(z),z_t)\,dt 
+ \alpha \int_0^T (\tfp(v),v_t)\,dt}_{\textsc{Forces III}}\,.
\end{multline*}
Substituting the above expression into \eqref{e:coupled-est-3} gives
\begin{multline} \label{e:almost-final} 
\big[t_1-t_0 +2C_{\tau,T,\epsilon_0}-2C_Te^{-\delta\tau}-2\eps\big] E(T)
+ \epsilon_0\int_0^T E(t)\, dt
\\
\le \big[2 C_{\tau,T,\epsilon_0}+ 2C_Te^{-\delta\tau}+2\epsilon_0 T+2C\big]E(0)+ lot(z,v)+
\\
\myspace + 2 \epsilon_0 \big\{\textsc{Forces I}\big\}+ \big\{\textsc{Forces II}\big\}
+  C_{\tau,T,\epsilon_0}\big\{\textsc{Forces III}\big\}\,.
\end{multline}

\medskip

\noindent 7. To prove that \eqref{e:almost-final} yields \eqref{e:pre-final},
we seek to render 
\begin{equation}\label{e:goal}
t_1-t_0 +2C_{\tau,T,\epsilon_0}-2C_Te^{-\delta\tau}-2\eps > 
2 C_{\tau,T,\epsilon_0}+ 2C_Te^{-\delta\tau}+2\epsilon_0T+2\eps\,,
\end{equation}
thus obtaining   
\begin{equation*}
\sigma:=\frac{2 C_{\tau,T,\epsilon_0}+ 2C_Te^{-\delta\tau}+2\epsilon_0T+2\eps }
{t_1-t_0 +2C_{\tau,T,\epsilon_0}-2C_Te^{-\delta\tau}-2 \eps}<1\,,
\end{equation*}
which readily simplifies to
\begin{equation}\label{e:simplified-goal}
t_1-t_0 -4C_Te^{-\delta\tau}> 2\eps +2\epsilon_0T\,.
\end{equation}
which holds for small $\eps,\eps_0$, and sufficiently large $\tau$.

Consequently, \eqref{e:almost-final} gives \eqref{e:pre-final}, where the interval 
$(t_0, t_1)$ has been translated to  a corresponding subset of $[s, s+T]$. 
Finally note that the result of Proposition \ref{p:pre-final} is continuous with respect to the finite-energy
topology, hence extends to all weak solutions. The proof of Proposition \ref{p:pre-final} is now complete.
\end{proof}

\bigskip

In order to move on from Proposition~\ref{p:pre-final} to conclusion that the hypothesis of Lemma~\ref{l:final} is satisfied for $T_{\gam_1,\gam_2}=\infty$ (independently of the trajectories) we must establish a bound on the products $(\tfw(z),z_t)$ and $(\tfp(v),v_t)$  appearing on the RHS of \eqref{e:pre-final}.

Such estimates can be derived using the structure of the  source terms $f_j$ and the energy identity (which, in particular, implies the integrability of the dissipation terms over $t\in (0,\infty)$). The   analysis of $(\tfp(v),v_t)$ can be carried out using this strategy (e.g. see \cite[Ch.~4 and 7]{chu-las-memoirs}). However, even though the plate is subject to full interior damping, in the wave component the norm $\|z_t\|^2$ cannot be controlled by $g(z_t)z_t$ alone, since the latter is restricted by the cutoff map 
$\chi(x)$ to just a subset of $\Om$.

Due to the geometric ``deficiency'' of the wave counterpart we follow the strategy employed in \cite{chu-las-tou}. The approach takes advantage of compactness of the attractor; this method was originally introduced in \cite{khan:b} to study von K\'arm\'an equation with \emph{internal} damping, and then later used for wave equation with boundary dissipation \cite{chu-las:snowbird} and von K\'arm\'an equations \cite{chu-las:jde:2007}. See also \cite{chu-las-memoirs} for an abstract realization of this idea.

\begin{proposition}\label{p:forces}
Suppose the assumptions of Proposition \ref{p:pre-final} hold. Then
the hypothesis of Lemma \ref{l:final} is satisfied with  $T_{\gam_1,\gam_2}=\infty$ for any two trajectories $\gam_1$, $\gam_2$ through the attractor.
\end{proposition}

\begin{proof}

The primary goal of this argument is to start with the estimate \eqref{e:pre-final} and 
get rid of the source-related terms collectively labeled $\{\textsc{Forces}\}$
in \eqref{d:forces}.

\medskip

\textbf{Step~1: Smoothness near stationary points.} The next inequality essentially invokes the  Lipschitz property of the derivative $f_1'$ of $f_1$ that follows from the second-order differentiability of $f_1$
\begin{equation}\label{p:forces:f1split}
\begin{split}
\intst (\tfw(z), z_t)
=&
(\tfw(z), z) \Big|_s^t - 
\intst\! \intOm z_t z^2 \Big(\int_0^1  \tfw''(\lam h + (1-\lam)\zet ) d\lam\Big) \, dxd\tau\,.
\end{split}
\end{equation}
The first term on the RHS is of a lower order, for the second term one can employ H\"{o}lder estimates, and  use the embedding  bound $\|z\|_6^2 \leq c E_{z,v}(t)$ to conclude:
\begin{equation}\label{p:forces:f1bound-A}
\Big|\intst (\tfw(z), z_t)\Big| \leq \eps( E_{z,v}(s)+E_{z,v}(t)) + C_{\eps}\sup_{\th\in[s,t]}\|z(\th)\|^2 + C_{\mathscr{A}} \intst \|z_t(\tau)\| E_{z,v}(\tau) d\tau\,.
\end{equation}
A similar decomposition can be repeated for $f_2$ which, however, requires a more general approach. We quote \cite[pp.~98--99]{chu-las-memoirs} to obtain
\[
\Big|\intst (\tfp(v), v_t) \Big|\leq 
C_{\mathscr{A}}\sup_{\th\in[s,t]}\|A_2^{1/2-\del} v\|^2 
+ C_{\mathscr{A}}  \intst \big( \|w_t(\tau)\| + \|v_t(\tau)\| \big)  E_{z,v}(\tau) \,d\tau\,.
\]
Further estimate  
\[
\big( \|w_t(\tau)\| + \|v_t(\tau)\| \big)  E_{z,v}(\tau) \leq 
C_{2,\mathscr{A},\eps} (\|w_t(\tau)\|^2 
+ \|v_t(\tau)\|^2)E_{z,v}(\tau) + \eps E_{z,v}(\tau)\,,
\] 
consequently the preceding estimate can be restated as 
\begin{equation}\label{p:forces:f2bound}
\Big|\intst (\tfp(v), v_t) \Big|\leq C_{\mathscr{A}} lot_s^t(z,v) 
+ C_{\mathscr{A},\eps} \! \intst \big( \|w_t\|^2 + \|v_t\|^2 \big) \,E_{z,v}(\tau) \,d\tau  
+ \eps \intst E_{z,v}(\tau) d\tau\,.
\end{equation}

\smallskip
Now pick a trajectory $\gam: t\mapsto \{\zet(t), w(t)\}$ through the attractor and for small $h$ define as before $\{z^h, v^h\} \dfn
\{\zet(t+h)-\zet(t), w(t+h)-w(t)\}$. From Theorem \ref{thm:attractor-existence} it follows that as $t\to \pm\infty$ this trajectory approaches the set of stationary points. 
In particular, the norm of the velocity components can be made as small as we wish, so let  $T_\gam^\eps$ be such that for all $h\in [0,1)$
\[
 \|\zet_t(t+h)\|,\; \|w_t(t+h)\| \leq 
\frac{\eps}{2(C_{\mathscr{A},\eps} +C_{\mathscr{A}})}
\quad \text{whenever}\quad t \leq T_{\gam}^\eps\,. 
\]
Substitute this estimate into  \eqref{p:forces:f1bound-A} and \eqref{p:forces:f2bound} (for $z=z^h$ and $v=v^h$) to conclude
\[
\begin{split}
\{\textsc{Forces}\}
\leq\; & \eps C_T \big[E_{z^h,v^h}(s)+E_{z^h,v^h}(t)\big] \\
&+ 2\eps C_T\intst E_{z^h, v^h}(\tau)\,d\tau +  C_{\mathscr{A},\eps,T} lot_s^t(z^h, v^h)\,.
\end{split}
\]
Now, let $\eps$ be small enough so that $\eps C_T  < \frac{1}{2}\min\left\{c_E, 1-\sig \right\}$ where the  constant $c_E$ comes from the result \eqref{e:pre-final} of  Proposition \ref{p:pre-final}. For $\sigma$ from $\eqref{e:pre-final}$ define
\[
\sig_1 \dfn \frac{\sig+C_T\eps}{(1-C_T\eps)}<1\quad \text{and}\quad 2C_T\eps < c_E\,.
\]
Combining the established estimate on $\{\textsc{Forces}\}$ with \eqref{e:pre-final} 
we obtain (for $t=s+T$)
\[
E_{z,v}(s+T)  \le \sig_1 E_{z,v}(s) 
+C_{\mathscr{A},\sigma,\eps,T}\, lot_{s}^{s+T}(z,v),\quad s < T_{\gam}^\eps - T\,.
\]
Since $\eps$ is now fixed we have obtained the hypothesis of  Lemma \ref{l:final} with finite $T_{\gam_1,\gam_2}$ and $G(\tau)\equiv 0$. Invoking this Lemma we can conclude that any trajectory through the attractor is strong, in the sense that
\begin{equation}\label{p:forces:smooth}
\{\zet,\zet_t, w, w_t\} \in H^2(\Om)\times H^1(\Om) \times H^4(\Gam_0)\times H^2(\Gam_0)\quad \forall t
\end{equation}
(however, possibly not globally bounded in this higher topology).
It remains to prove that the RHS terms of \eqref{p:forces:f1bound-A} and \eqref{p:forces:f2bound} can be estimated as in the hypothesis of  Lemma \ref{l:final}, but now for $T_{\gam_1,\gam_2}=\infty$.

\medskip

\textbf{Step 2: Analysis of the source for the plate component \eqref{p:forces:f2bound}.} Since the plate is damped on all of its interior it is possible to apply a standard argument \cite[Section 4.1]{chu-las-memoirs} that provides control over 
$\int_s^t (\|w_t\|^2+\|u_t\|^2)$ via the dissipation. From the Assumption \ref{a:dampings-1} we have
\[
m'_b\intst \|w_t\|^2 \leq  \intst \int_{\Gam_0} b(w_t)w_t\,.
\]
The global bounds on the energy (Proposition \ref{p:energy-bound}) and the energy identity \eqref{e:energy-identity} verify that $\int_s^t \|w_t\|^2$ is bounded uniformly for 
$s, t\in \bbR$ (for trajectories through the attractor). 
Consequently, in \eqref{p:forces:f2bound} we may define
\begin{equation}\label{p:forces:d}
G(\tau) \dfn  C_{\scrA,\eps}\, \big(\|w_t(\tau)\|^2 + \|v_t(\tau)\|^2\big)
\end{equation}
which belongs to $L^1(\mathbb{R})$.

\medskip

\textbf{Step 3: Analysis of the source for the wave component \eqref{p:forces:f1bound-A}}
The argument employed for the plate cannot be repeated here because the acoustic damping is only supported on a subset of $\Om$. However now we can use the fact that the trajectories through the attractor are smooth \eqref{p:forces:smooth}, and that the attractor itself $\mathscr{A}$ is compact in the finite energy space $Y$.

From the compactness property we know that given any $\eps>0$ there exists a \emph{finite} set $\{\zet^j_\eps, u^j_\eps, w^j_\eps, y^j_\eps\}_{j=1}^N\subset \mathscr{A}$ such that for any fixed time $t$ one can find $i=i(t)$ and $j=j(t)\in \bbN$, $1\leq i,j\leq N$ for which the difference $u^{i,j}_\eps\dfn u^i_\eps-u^j_\eps$ satisfies
\[
\|z_t(t) - u^{i(t),j(t)}_\eps\| < \eps\,.
\]
Moreover, since the collection $u_\eps^j$ is finite and belongs to the attractor, then
\[
\sup_{i,j} \|u^{i,j}_\eps\|_{1,\Om} \leq C_{\mathscr{A},\eps}\,.
\]
It is now possible to refine the decomposition \eqref{p:forces:f1bound-A} of the acoustic source: for any fixed time $t=\tau$, we have
\begin{equation}\label{p:forces:f1bound-b1}
\begin{split}
&\intOm z_t z^2 \Big(\int_0^1  \tfw''(\lam h + (1-\lam)\zet ) d\lam\Big) \, dx\\
=&\intOm (z_t-u^{i,j}_\eps + u^{i,j}_\eps) z^2 \Big(\int_0^1  \tfw''(\lam h + (1-\lam)\zet ) d\lam\Big) \, dx\\
\leq &
\intOm \bigg(|z_t-u^{i,j}_\eps| + |u^{i,j}_\eps|\bigg) z^2 \Big|\int_0^1  \tfw''(\lam h + (1-\lam)\zet ) d\lam\Big| \, dx\,.
\end{split}
\end{equation}
Since $z_t(\tau)$ belongs to the attractor, then it is possible to choose  $i=i(\tau), j=j(\tau)$ such that $u_\eps^{i,j}$ satisfies
\begin{equation}\label{p:forces:f1bound-b2}
\begin{split}
\intOm |z_t-u^{i,j}_\eps| z^2 \left|\int_0^1  \tfw''(\lam h + (1-\lam)\zet ) d\lam\right| \, dxd\tau
&\leq C_{\mathscr{A}}\, \|z_t(\tau)-u^{i(\tau),j(\tau)}_\eps\| E_{z,v}(\tau)\\
&\leq C_{\mathscr{A}}\, \eps E_{z,v}(\tau)\,,
\end{split}
\end{equation}
whereas when applying H\"{o}lder estimates to the remaining terms in \eqref{p:forces:f1bound-b1} we take advantage of the fact that $u^{i,j}_\eps \in H^1(\Om)$, the norm being bounded over all (finitely many) 
$i,j$:
\[
\intOm |u^{i,j}_\eps| z^2 \Big|\int_0^1  \tfw''(\lam h + (1-\lam)\zet ) d\lam\Big| \, dxd\tau
\leq
C\|u^{i,j}_\eps\|_p  \|z\|_{6-\del}^2 \big\|1+ |h|+|\zet|\big\|_{\frac{(6-\del)p'}{6-\del-2p'}}
\]
for conjugate exponents $p$ and $p'$.  
Pick $\del>0$ and $p >2$  so that $\frac{(6-\del)p'}{6-\del-2p'} \leq  6$, for instance: $p=3$ (whence $p'=3/2$) and $\del =1$. 
We note that  
\begin{itemize}
		\item since $\frac{(6-\del)p'}{6-\del-2 p'} =  15/4\leq 6$ we may bound the corresponding norm  by $C_{\mathscr{A}}$;

	\item  $\|u^{i,j}_\eps\|_{p} \leq \|u^{i,j}\|_{1,\Om} \leq C_{\scrA,\eps}$ for all $i,j$;

	\item  using Sobolev embeddings and interpolation gives
	\begin{equation}\label{p:forces:f1bound-b3}
	\|z\|_{6-\del}^2\leq \frac{\eps}{C\cdot C_{\mathscr{A}} \cdot 
C_{\mathscr{A},\eps}} E_{z,v}(\tau) + C_{2,\mathscr{A},\eps} \|z(\tau)\|^2\,.
	\end{equation}
\end{itemize}
Combining these results, we conclude (relabeling the constants)
\[
\intOm |u^{i,j}_\eps| z^2 \Big|\int_0^1  \tfw''(\lam h + (1-\lam)\zet ) d\lam\Big|\,dxd\tau
\leq C_{3,\mathscr{A},\eps} \|z(\tau)\|^2 + \eps C_{\scrA} E_{z,v}(\tau)\,.
\]
So, the decomposition \eqref{p:forces:f1split} refined with \eqref{p:forces:f1bound-b1}, \eqref{p:forces:f1bound-b2}, \eqref{p:forces:f1bound-b3} yields
\begin{equation}\label{p:forces:f1bound-B} 
\left|\intst (\tfw(z), z_t)\right| \leq \eps C_{\mathscr{A}}( E_{z,v}(s)+E_{z,v}(t)) + C_{\eps,\mathscr{A}}\sup_{\th\in[s,t]}\|z(\th)\|^2 + \eps \intst  E_{z,v}(\tau)d\tau\,.
\end{equation}
Thus, the extra regularity of the attractor combined with its compactness in the finite energy space show that, despite the geometrically restricted dissipation, the interaction of the acoustic source with pressure  $z_t$ is (``almost'') of a lower order.

\smallskip
At this stage the last estimate \eqref{p:forces:f1bound-B} on the source $f_1$, along with the decomposition of the structural source  \eqref{p:forces:f2bound}, and the integrability of the velocity (controlled by the full-interior dissipation)  \eqref{p:forces:d} can be substituted into  the perturbed observability estimate \eqref{e:pre-final}. For small enough $\eps$ we obtain the hypothesis \eqref{l:final:assum} of  Lemma \ref{l:final} with no restriction on $T_{\gam_1,\gam_2}$.  The proof of  Proposition \ref{p:forces} is now complete.
\end{proof}

\begin{remark}
	\begin{rm}
	It is also possible to employ the approach used for the acoustic source $f_1$ (the final  step in the proof of Proposition \ref{p:forces}) to derive an analogous estimate for the plate $f_2$, instead of appealing to the integrability of $\intst (b(w_t),w_t)$ (as was done in Step 2 of that proof). For that one can go back to the derivation of  \eqref{p:forces:f2bound} (see \cite[pp. 98--99]{chu-las-memoirs}) and rewrite $v_t = v_t - y_\eps^{ij}+y_{\eps}^{ij}$ where $y_{\eps}^{ij} = y^i_\eps - y^j_{\eps}$ is a suitable set of points approximating the compact attractor $\scrA$.
	\end{rm}
\end{remark}

\medskip
Via Proposition \ref{p:forces}, the part (a) of  Lemma \ref{l:final} now implies  that the observability estimate \eqref{e:final} holds for all $s\in\mathbb{R}$ independently of the chosen trajectories through the attractor, which is the last ingredient necessary to establish finite-dimensionality of $\mathscr{A}$. Thus, the part (b) of Lemma \ref{l:final} now shows that the bound on the higher energy of the trajectories is global and independent of the trajectories themselves, confirming that $\mathscr{A}$ is a bounded subset of the strong topology. This step completes the proof of Theorem \ref{thm:finite-dim}.

\vskip 1cm
\appendix

\section{Basic definitions and major abstract results}
\label{appendix-a}
We begin by recording two general (abstract) theorems which will be invoked 
in the proof of our main results. The first theorem describes a sufficient condition 
to ensure asymptotic smoothness of a semi-flow; this property is crucial 
for the {\em existence} of a global attractor.
The second result plays a major role in showing that the obtained attractor 
has {\em finite fractal dimension}.
 
\begin{definition}\label{d:asymptotic}
A dynamical system $(X,S(t))$ is said to be {\em asymptotically smooth}
if for any bounded set $B$ which is forward invariant (i.e. $S(t)B\subset B$, 
$t\ge 0$) there exists a compact set $K\subset\overline{B}$ such that
\begin{equation*}
{\rm dist}(S(t)B,K) \rightarrow 0\,, \qquad 
\textrm{as} \quad t\to +\infty\,.
\end{equation*}
\end{definition}

\begin{theorem}[Proposition~2.10 in \cite{chu-las-memoirs}]
\label{t:khan}
Let $(X,S(t))$ be a dynamical system on a complete metric space $X$, with metric $d$.
Assume that for any bounded positively invariant set $B\subset X$ and any $\epsilon>0$
there exists $T=T(\epsilon,B)$ such that 
\begin{equation*} 
d(S_TY_1,S_TY_2)\le \epsilon + \Psi_{\epsilon,B,T}(Y_1,Y_2)\,,
\qquad y_i\in B\,,
\end{equation*}
where $\Psi_{\epsilon,B,T}(Y_1,Y_2)$ is a nonnegative function defined on $B\times B$ 
such that
\begin{equation} \label{e:sequential-limit}
\liminf_{m\to \infty} \liminf_{n\to \infty}\Psi_{\epsilon,B,T}(y_n,y_m)=0
\end{equation}
for every sequence $\{y_n\}_n$ in $B$.
Then the dynamical system $(X,S_t)$ is asymptotically smooth.
\end{theorem}

Let us recall the classical definition of {\em fractal} (or box-counting) dimension 
of a compact set; see \cite{temam-97}, or \cite{babin}.

\begin{definition}[Fractal dimension]
The fractal dimension $\dim_f E$ of a compact set $E$ is defined by
\begin{equation*}
\dim_f E:=\limsup_{r\to 0^+}\frac{\ln n(E,r)}{\ln (1/r)}\,,
\end{equation*}
where $n(E,r)$ is the minimal number of closed sets of diameter $2r$ which cover $E$.
\end{definition}

The following result provides a generalization of a Ladyzhenskaia's theorem
on dimension of invariant sets; its requirements are somehow related to the so 
called ``smoothing-squeezing'' property introduced in \cite{prazak}.
It should also be noted that the forthcoming estimate of the fractal dimension given by 
\eqref{e:estimate-dimension} is not sharp.

For general criteria yielding effective estimates of the fractal dimension of a global
attractor, as well as the application of these results to various specific equations,
see \cite{temam-97}.
More recent advances include, e.g., \cite{chepyzhov-ilyin}, yielding a sharper estimate 
of the dimension of the attractor for the 2D Navier-Stokes equations; see also \cite{constantin-levant-titi:06}.
It should be observed however that this kind of results is inapplicable in the present context, because of the lack of differentiability of the flow.

\begin{theorem}[Theorem~2.15 in \cite{chu-las-memoirs}]
\label{t:ladyzhenskaia-general}
Let $X$ be a Banach space and $M$ be a bounded closed set in $X$.
Assume that there exists a mapping $V:M \to X$ such that 
\begin{enumerate} 
\item[(i)]
$M\subseteq V M$;
\item[(ii)]
$V$ is Lipschitz on $M$, i.e. there exists $L>0$ such that 
\begin{equation*}
\|Vv_1-Vv_2\|\le L \|v_1-v_2\|\,, \quad v_1,v_2\in M\,;
\end{equation*}
\item[(iii)]
there exist compact seminorms $n_1(\cdot)$ and $n_2(\cdot)$ on $X$ satisfying 
\begin{equation}\label{e:involves-seminorms}
\|Vv_1-Vv_2\|\le \eta \|v_1-v_2\|+K[n_1(v_1-v_2)+n_2(Vv_1-Vv_2)]
\end{equation}
for any $v_1, v_2\in M$, where $0<\eta<1$ and $K>0$ are constants.
(A seminorm $n(\cdot)$ on $X$ is said to be compact if for any 
bounded set $B\subset X$ there exists a sequence $\{x_n\}\subset B$ such that 
$n(x_m-x_n)\to 0$, as $n,m\to \infty$.)
Then $M$ is a compact set in $X$ with finite fractal dimension.
\end{enumerate}
Moreover, one has the estimate 
\begin{equation}\label{e:estimate-dimension}
dim_f M\le \log m_0\Big(\frac{4K(1+L^2)^{1/2}}{1-\eta}\Big)
\Big(\log \frac2{1+\eta}\Big)^{-1}\,,
\end{equation}
where $m_0(R)$ is the maximal number of pairs $(x_i,y_i)$ in $X\times X$
possessing the properties 
\begin{equation*}
\|x_i\|^2+\|y_i\|^2\le R^2\,, \quad n_1(x_i-x_j)+n_1(y_i-y_j)>1\,, \; i\ne j\,.
\end{equation*}
\end{theorem}

\section{Accessory proofs}

\smallskip
\noindent
{\em Proof of the estimate \eqref{e:final-for-challenging}. }
Let us recall the definition of
\begin{equation*}
I_1:=\int_{Q_T}\psi z\tfw(z)\vh\cdot \nabla \psi\,.
\end{equation*}
Since readily
\begin{equation*}
\big|I_1\big|\le C\int_0^T \int_{\Omega}|z\tfw(z)|\,,
\end{equation*}
to show the inequality \eqref{e:final-for-challenging} 
we need to estimate the integral $J(t)$ ($J$, in short) defined by
\begin{equation}
J:= \int_{\Omega}|z\tfw(z)|= \int_{\Omega}|z|\,|f_1(\zeta+z)-f_1(z)|
= \int_{\Omega}|z|\,|f(\zeta+z)-f(z)+\lambda z|\,.
\end{equation}
By using first Assumption~\ref{h:source-1} on the nonlinear function $f$ 
(which occurrs in the wave equation), in the form
\begin{equation*}
|f(s_1)-f(s_2)|\le c \,(1+ |s_1|^q+  |s_2|^q)\, |s_1-s_2| 
\qquad \forall s_1, \, s_2\; (q\le 2)\,,
\end{equation*}
and applying next the H\"older inequality with exponents $(3,3/2)$, we obtain
\begin{align}
J &\le c_1\,\int_{\Omega}\big(1+|\zeta+z|^q+|\zeta|^q\big)\,|z|^2
\le c_1\Big[\int_{\Omega}\big(1+|\zeta+z|^q+|\zeta|^q\big)^3\Big]^{1/3}\,
\left[\int_{\Omega}|z|^3\right]^{2/3}
\nonumber\\[1mm]
& \le c_2 \,\Big[\int_{\Omega}\big(1+|\zeta+z|^{3q}+|\zeta|^{3q}\big)\Big]^{1/3}\,
\|z\|_{L_3(\Omega)}^2
\nonumber\\[1mm]
& \le c_3 \,\big(1+\|\zeta+z\|_{L_{3q}(\Omega)}^q+\|\zeta\|_{L_{3q}(\Omega)}^q\big)\,
\|z\|_{L_3(\Omega)}^2
\label{e:trivial-for-everybody}
\end{align}
Note that to get the last inequality we used that 
$\|\zeta\|^q \le 1+ \|\zeta\|^2$ for any $\zeta$, since $q\le 2$.
Then, using the Sobolev embeddings $H^1(\Omega)\subset L_6(\Omega)\subset L_{3q}(\Omega)$
and $H^{1/2}(\Omega)\subset L_3(\Omega)$, we see that 
\eqref{e:trivial-for-everybody} implies
\begin{equation*}
J\le c_3 \big(1+\|\zeta+z\|_{H^1(\Omega)}^2+\|\zeta\|_{H^1(\Omega)}^2\big) C_{\Om}
\|z\|_{H^{1/2}(\Omega)}^2\,.
\end{equation*}
Thus, since $\|\zeta+z\|_{1,\Omega}$ and $\|\zeta\|_{1,\Omega}$
are uniformly (in $t$) bounded for initial data in a bounded set $\sB$,
we finally obtain 
\begin{equation*}
J(t)\le c_3 (3+C_{\sB})\,C_{\Om}\|z\|_{H^{1/2}(\Omega)}^2\,,
\end{equation*}
which combined with the interpolation inequality 
\begin{equation*}
\|z\|_{H^{1/2}(\Omega)}\le 
\epsilon_0\|z\|_{H^{1}(\Omega)}+ C_{\epsilon_0}\|z\|_{L_2(\Omega)}
\end{equation*}
finally yields the estimate \eqref{e:final-for-challenging}.
\qed


\end{document}